\documentclass[a4paper]{amsart} % sets paper format is a4
\usepackage[ascii]{inputenc} % for compatibility, please do not use utf8, let alone latin1
\usepackage{amssymb}
\usepackage{mathrsfs}
\usepackage[hidelinks]{hyperref}

\usepackage{xcolor}

\newcommand{\red}[1]{\textcolor{red}{#1}}

\newcommand{\cyan}[1]{\textcolor{cyan}{#1}}

\usepackage[shortlabels]{enumitem}
\setlist[enumerate,1]{label={(\Alph*)}}
\setlist[enumerate,2]{label={(\alph*)}}
\setlist[enumerate,3]{label={$\bullet_{\arabic*}$}}

\newenvironment{PROOF}[2][\proofname.]
   {\begin{proof}[#1]\renewcommand{\qedsymbol}{\textsquare$_{\mathrm{#2}}$}}
   {\end{proof}}

\newtheorem{theorem}{Theorem}[section] 
\newtheorem{claim}[theorem]{Claim}
\newtheorem{cd}[theorem]{Claim/Definition}

\newtheorem{observation}[theorem]{Observation}

\theoremstyle{definition}
\newtheorem{definition}[theorem]{Definition}
\newtheorem{dc}[theorem]{Definition/Claim}

\newtheorem{discussion}[theorem]{Discussion}

\newtheorem{hypothesis}[theorem]{Hypothesis}

\theoremstyle{remark}
\newtheorem{remark}[theorem]{Remark}

\newtheorem{notation}[theorem]{Notation}
\newtheorem{context}[theorem]{Context}

\newcommand{\acl}{\mathrm{acl}}
\newcommand{\dcl}{\mathrm{dcl}}

\newcommand{\Ax}{\mathrm{Ax}}
\newcommand{\bas}{\mathrm{bas}}

\newcommand{\brim}{\mathrm{brim}}

\newcommand{\Cb}{\mathrm{Cb}}

\newcommand{\ess}{\mathrm{ess}}
\newcommand{\full}{\mathrm{full}}
\newcommand{\good}{\mathrm{good}}
\newcommand{\na}{\mathrm{na}}

\newcommand{\ortp}{\mathrm{ortp}}
\newcommand{\PC}{\mathrm{PC}}

\newcommand{\qr}{\mathrm{qr}}

\newcommand{\reg}{\mathrm{reg}}

\newcommand{\stp}{\mathrm{stp}}

\newcommand{\vq}{\mathrm{vq}}

\newcommand{\wg}{\mathrm{wg}}
\newcommand{\wk}{\mathrm{wk}}

\newcommand{\oor}{\mathrm{or}}

\newcommand{\genw}{\mathrm{gw}}

\newcommand{\arity}{\mathrm{arity}}

\newcommand{\cf}{\mathrm{cf}}

\newcommand{\dom}{\mathrm{dom}}

\newcommand{\id}{\mathrm{id}}

\newcommand{\inv}{\mathrm{inv}}

\newcommand{\Rang}{\mathrm{Rang}}

\newcommand{\seq}{\mathrm{seq}}

\newcommand{\bs}{\mathrm{bs}}

\newcommand{\EM}{\mathrm{EM}}
\newcommand{\eq}{\mathrm{eq}}

\newcommand{\fin}{\mathrm{fin}}

\newcommand{\LST}{\mathrm{LST}}

\newcommand{\pr}{\mathrm{pr}}

\newcommand{\tp}{\mathrm{tp}}

\newcommand{\bfa}{\mathbf{a}}

\newcommand{\bfe}{\mathbf{e}}
\newcommand{\bfF}{\mathbf{F}}

\newcommand{\bfi}{\mathbf{i}}
\newcommand{\bfJ}{\mathbf{J}}

\newcommand{\bfk}{\mathbf{k}}

\newcommand{\bfP}{\mathbf{P}}

\newcommand{\bbE}{\mathbb{E}}   % note:  not \bE

\newcommand{\mn}{\medskip\noindent}
\newcommand{\sn}{\smallskip\noindent}

\newcommand{\cE}{\mathscr{E}}

\newcommand{\cP}{\mathscr{P}}

\newcommand{\clS}{\mathcal{S}}

\newcommand{\gC}{\mathfrak{C}}

\newcommand{\gK}{\mathfrak{K}}

\newcommand{\gk}{\mathfrak{k}}

\newcommand{\gn}{\mathfrak{n}}

\newcommand{\gs}{\mathfrak{s}}
\newcommand{\gt}{\mathfrak{t}}

\newcommand{\eps}{\varepsilon}
\newcommand{\lh}{{\ell g}}
\newcommand{\rest}{\restriction}

\newcommand{\caret}{{\char 94}}
\newcommand{\LL}{\langle}
\newcommand{\RR}{\rangle}

\newcommand{\olsi}[1]{\,\overline{\!{#1}}} % overline short italic

\newcommand*{\defeq}{\mathrel{\vcenter{\baselineskip0.5ex \lineskiplimit0pt\hbox{\scriptsize.}\hbox{\scriptsize.}}}=}

\newcount\skewfactor
\def\mathunderaccent#1#2 {\let\theaccent#1\skewfactor#2
\mathpalette\putaccentunder}
\def\putaccentunder#1#2{\oalign{$#1#2$\crcr\hidewidth
\vbox to.2ex{\hbox{$#1\skew\skewfactor\theaccent{}$}\vss}\hidewidth}}

\newbox\noforkbox \newdimen\forklinewidth
\setbox0\hbox{$\textstyle\bigcup$}
\setbox1\hbox to \wd0{\hfil\vrule width \forklinewidth depth \dp0
                        height \ht0 \hfil}
\setbox\noforkbox\hbox{\box1\box0\relax}
\def\unionstick{\mathop{\copy\noforkbox}\limits}
\def\nonfork#1#2_#3{#1\unionstick_{\textstyle #3}#2}
\def\nonforkin#1#2_#3^#4{#1\unionstick_{\textstyle #3}^{\textstyle
    #4}#2}
\setbox0\hbox{$\textstyle\bigcup$}
\setbox1\hbox to \wd0{\hfil{\sl /\/}\hfil}
\setbox2\hbox to \wd0{\hfil\vrule height \ht0 depth \dp0 width
                                \forklinewidth\hfil}
\newbox\doesforkbox
\setbox\doesforkbox\hbox{\box1\box0\relax}
\def\nunionstick{\mathop{\copy\doesforkbox}\limits}

\def\fork#1#2_#3{#1\nunionstick_{\textstyle #3}#2}
\def\forkin#1#2_#3^#4{#1\nunionstick_{\textstyle #3}^{\textstyle
    #4}#2}

\newcommand{\stickT}{%
\setbox255=\hbox{\raise1ex\hbox{$\hspace{0.2pt}\,\bullet\,$}}
\mathord{\rlap{\hbox to\wd255{\hss\hbox{$|$}\hss}}
\box255}
}
\newcommand{\stickS}{%
\setbox255=\hbox{\raise0.6ex\hbox{$\scriptstyle\bullet$}}
\mathord{\rlap{\hbox to\wd255{\hss\hbox{$\scriptstyle|$}\hss}}
\box255}
}

\begin{document}
\makeatletter\def\shfiuwefootnote{\gdef\@thefnmark{}\@footnotetext}\makeatother\shfiuwefootnote{Version 2023-03-24. See \url{https://shelah.logic.at/papers/839/} for possible updates.}

\title {Stable frames and weights \\
 Sh839}
\author {Saharon Shelah}
\address{Einstein Institute of Mathematics\\
Edmond J. Safra Campus, Givat Ram\\
The Hebrew University of Jerusalem\\
Jerusalem, 91904, Israel\\
 and \\
 Department of Mathematics\\
 Hill Center - Busch Campus \\ 
 Rutgers, The State University of New Jersey \\
 110 Frelinghuysen Road \\
 Piscataway, NJ 08854-8019 USA}
\email{shelah@math.huji.ac.il}
\urladdr{http://shelah.logic.at}
\thanks{For earlier versions, the author thanks Alice Leonhardt for the beautiful typing up to 2019.
The author would like to thank the Israel Science Foundation for
partial support of this research (Grant No. 242/03). 
First Typed 2002-Sept.-10.\\
For later versions, the author would like to thank the NSF and BSF for partially supporting this research --- NSF-BSF 2021: grant with M. Malliaris, NSF 2051825, BSF 3013005232 (2021/10-2026/09). The author is also grateful to an individual who wishes to remain anonymous for generously funding typing services, and thanks Matt Grimes for the careful and beautiful typing.\\
References like [Sh:950, Th0.2=Ly5] mean that the internal label of Th0.2 is y5 in Sh:950.
The reader should note that the version in my website is usually more up-to-date than the one in arXiv.}

% Previous version - 2012/May/22

% formerly F569 - changed March 2004
% formerly titled ``More on Frames"

% The following is necessary as the current amsart.cls does not allow 2020
% with amsart.cls starting version 2.20.6 (contained in TeX Live 2020) 
% the following definition should be removed
\makeatletter
\@namedef{subjclassname@2020}{\textup{2020} Mathematics Subject Classification}
\makeatother


\date{March 7, 2023}

\begin{abstract}

\noindent
We would like to generalize imaginary elements, 
weight of $\ortp(a,M,N)$, $\bfP$-weight, $\bfP$-simple types, etc. from
\cite[Ch.III,V,\S4]{Sh:c} to the context of good frames.  This requires
allowing the vocabulary to have predicates and function symbols of
infinite arity, but it seems that we do not suffer any real loss.
\end{abstract}

\maketitle
\numberwithin{equation}{section}
\setcounter{section}{-1}

\newpage

\centerline {Annotated Content}
\bigskip

%\noindent
%Part I, \hfill pg.\pageref{part 1}
%\bigskip

\noindent
\S0 \quad Introduction \hfill pg.\pageref{0}
%\mn
%\begin{enumerate}
%\item[]  [   ]
%\end{enumerate}
%\bigskip
\bigskip

\noindent
\S1 \quad Weight and $\bfP$-weight  \hfill pg.\pageref{1}
\mn
\begin{enumerate}
\item[]  [For $\gs$ a good $\lambda$-frame with some additional
  properties, we define placed and $\bfP$-weight.]
\end{enumerate}
\bigskip

\noindent
\S2 \quad Imaginary elements, an $\ess$-$(\mu,\lambda)$-AEC and
frames  \hfill pg.\pageref{2}
\mn
\begin{enumerate}
\item[]  [Define an $\ess$-$(\mu,\lambda)$-AEC allowing
infinitary functions.  Then get $\gs$ with type bases.]
\end{enumerate}
\bigskip

\noindent
\S3 \quad $\bfP$-simple types 
\hfill pg.\pageref{3}

%\noindent Private Appendix, pg.\pageref{PA} % removed by martin 2020-04-17
\newpage

%\centerline {Part I: Beautiful frames: weight and simplicity} \label{part 1}

\section {Introduction} \label{0}

%We consider here the directions listed in the abstract\footnote{As
%we have started this in 2002 and have not worked on it for long, we
%intend to make public what is in reasonable state.}

We assume $\gs$ is a good $\lambda$-frame with some
extra properties from \cite{Sh:705} (e.g., as in the assumption of
\cite[\S12]{Sh:705}) so we shall assume knowledge of \cite{Sh:705} and
the basic facts on good $\lambda$-frames from \cite{Sh:600}. 

We can look at results from \cite{Sh:c} which were 
not regained in beautiful $\lambda$-frames.  Well, of
course, we are far from the main gap for the original $\gs$ 
(\cite[Ch.XIII]{Sh:c}) and there are
results which are obviously more strongly connected to elementary
classes, particularly ultraproducts.  This leaves us with parts of
type theory: regular and semi-regular types, weight, 
    $\bfP$-simple\footnote{The motivation is that for suitable $\bfP$ 
    (e.g. a single regular type),  on the one hand 
    $$\stp(a,A) \not\perp \bfP \Rightarrow ``\stp(a/E,A) \text{ is $\bfP$-simple for some
    equivalence relation definable over $A$"}$$ and on the other hand, 
    if $\stp(a_i,A)$ is $\bfP$-simple for $i < \alpha$ then 
    $\Sigma\big\{w(a_i,A) \cup \{a_j : j < i\}\red{)} : i < \alpha \big\}$
    does not depend on the order in which we list the $a_i$-s.  Note that 
    $\bfP$ here is ${\cP}$ there.} 
types, ``hereditarily orthogonal
to $\bfP$" (the last two were defined and investigated in 
\cite[Ch.V,\S0 + Def4.4-Ex4.15]{Sh:a},
\cite[Ch.V,\S0,pg.226,Def4.4-Ex4.15,pg.277-284]{Sh:c}). 

Some of Hrushovski's profound works are a continuation of \cite[\S4]{Sh:a} and \cite[\S4]{Sh:c}, but note that ``a type $q$ is $p$-simple (or $\bfP$-simple)" and ``$q$
is hereditarily orthogonal to $p$ (or $\bfP$)" here are essentially
    the\footnote{Note, ``foreign to $\bfP$" and ``hereditarily
    orthogonal to $\bfP$" are equivalent.  Now (with $\bfP = \{p\}$ for simplicity)
    \mn
    \begin{enumerate}[(a)]
        \item  $q(x)$ is $p(x)$-simple when for some set $A$\red{,} in $\gC$ we have 
        $q(\gC) \subseteq \acl(A \cup \bigcup p_i(\gC))$.
    \sn
        \item  $q(x)$ is $p(x)$-internal when for some set $A$, in $\gC$ we have $q(\gC) \subseteq \dcl(A \cup p(\gC))$.
    \end{enumerate}
    \mn
    Note
    \mn
    \begin{enumerate}
        \item[$(\alpha)$]   Internal implies simple.
    \sn
        \item[$(\beta)$]  If we aim at computing weights, it is better to stress acl as it covers more.
    \sn 
        \item[$(\gamma)$]   But the difference is minor, and in existence it is better to stress dcl.
    \sn
        \item[$(\delta)$]  Also, it is useful that 
        $$\big\{F \rest (p(\gC) \cup q(\gC) : F \text{ an automorphism of $\gC$ over } p(\gC) \cup \dom(p) \big\}$$
        is trivial when $q(x)$ is $p$-internal but not so for $p$-simple
        (though form a pro-finite group).
    \end{enumerate}} 
``internal" and ``foreign" there. %in Hrushovshi's profound works.

 For more on understanding regular types in the first order case, see both \cite{Sh:401} and Laskowski and the author in \cite{Sh:933}.

\bigskip
\begin{center}
    * \quad * \quad *
\end{center}

\bigskip
This paper was Part I of the original \cite{Sh:839}, which has existed (and circulated to some extent) since 2002. The second and third parts have been split off into \cite{Sh:1238}, \cite{Sh:1239}. They have been continued in \cite{Sh:851} and \cite{Sh:871}, respectively.

\begin{notation}
As in \cite{Sh:c}, \cite{Sh:e}, $M$ and $N$ are models, $M$ has vocabulary $\tau_M$, $|M|$ is its universe and $\|M\|$ its cardinality. We write $\ortp(-)$ for the orbital type.
\end{notation}

\newpage

\section {Weight and $\bfP$-weight} \label{1}

%Recalling \red{[4A]}%\cite{Sh:600}, \cite{Sh:705}:
On `good$^+$,' see Definition \cite[1.3(1), pg.7]{Sh:705} and Claim \cite[1.5(1), pg.7]{Sh:705}, which relies on \cite[\S3]{Sh:600}, \cite{Sh:576}.

\sn
On `type-full,' see Definition \cite[6.35, pg.112]{Sh:600}: it means $\clS_\gs(M) = \clS_\bfk^1(M)$.

\sn
On primes and $K_\gs^{3,\qr}$, see \cite[5.15, pg.73]{Sh:705}.

\sn
On $K_\gs^{3,\vq} = K_\gs^{3,\qr}$, see Definition \cite[5.9, pg.69]{Sh:705}.

\sn
On orthogonality, see \cite[\S6]{Sh:705}.

\begin{context}\label{a2}
1) $\gs$ is a type-full $\good^+ \,\lambda$-frame
with primes, $K^{3,\vq}_\gs = K^{3,\qr}_\gs$, ${\perp} = 
{\underset \wk \perp}$ and $p \perp M \Leftrightarrow
p\ {\underset {\mathrm{a}} \perp}\ M$. Note that as $\gs$ is full,
$\clS^\bs_\gs(M) = \clS^\na_\gs(M)$; also, $\gk_\gs =
\gk[\gs] = (K^\gs,\le_{\gk_\gs})$ is the AEC.

\noindent
2) $\gC$ is an $\gs$-monster so it is 
$K^\gs_{\lambda^+}$-saturated over $\lambda$, and $M <_\gs \gC$
means $M \le_{\gk[\gs]} \gC$ and $M \in K_\gs$.  
As $\gs$ is full, it has regulars.
\end{context}

\begin{observation}\label{a5}
$\gs^\reg$ satisfies all the above except being full.
\end{observation}

\begin{remark}\label{a6}
Recall $\gs^\reg$ is derived from $\gs$, replacing $\clS_\gs(M)$ by $\{p \in \clS_\gs(M) : p \text{ regular}\}$ (see \cite[10.18, pg.164]{Sh:705}).    
\end{remark}

\begin{proof}  See \cite[10.18=L10.p19tex]{Sh:705} and Definition
\cite[10.17=L10.p18tex]{Sh:705}. 
\end{proof}

\begin{claim}\label{a8}
1) If $p \in \clS^\bs_\gs(M)$ 
\underline{then} we can find $b,N$ and a finite $\bfJ$ such that:
\mn
\begin{enumerate}
    \item[$\circledast$]  
    \begin{enumerate}[$(a)$]
        \item $M \le_\gs N$

        \item $\bfJ \subseteq N$ is a finite independent set in $(M,N)$.

        \item $c \in \bfJ \Rightarrow \ortp(c,M,N)$ is regular (recalling that $\ortp$ stands for `orbital type').

        \item $(M,N,\bfJ) \in K^{3,\qr}_\gs$

        \item $b \in N$ realizes $p$.
    \end{enumerate} 
\end{enumerate}
\mn
2) If $M$ is brimmed, we can add
\mn
\begin{enumerate}
    \item[$(f)$]   $(M,N,b) \in K^{3,\pr}_\gs$.
\end{enumerate}
\mn
3) In (2), $|\bfJ|$ depends only on $(p,M)$.

\noindent
4) If $M$ is brimmed, \underline{then} we can work in $\gs(\brim)$ and get the
same $\|\bfJ\|$ and $N$ (so $N \in K_\gs$ {is} brimmed).
\end{claim}

\begin{PROOF}{\ref{a8}}  
1) By induction on $\ell < \omega$, we try to choose 
$N_\ell,a_\ell,q_\ell$ such that:
\mn
\begin{enumerate}
    \item[$(*)$]   
    \begin{enumerate}[$(a)$]
        \item $ N_0 = M$

        \item $ N_\ell \le_\gs N_{\ell +1}$

        \item $q_\ell \in \clS_\gs(N_\ell)$, so possibly 
        $q_\ell \notin \clS^\na_\gs(N_\ell)$.

        \item $q_0 = p$

        \item $q_{\ell +1} \rest N_\ell = q_\ell$

        \item $q_{\ell +1}$ forks over $N_\ell$, so now \underline{necessarily} 
        $q_\ell \notin \clS^\na_\gs(N_\ell)$.

        \item $(N_\ell,N_{\ell +1},a_\ell) \in K^{3,\pr}_\gs$

        \item $r_\ell = \ortp(a_\ell,N_\ell,N_{\ell +1})$ is regular.

        \item $r_\ell$  either is $\perp M$ or does not fork over $M$.
    \end{enumerate}
\end{enumerate}
\mn
If we succeed to carry the induction for all $\ell < \omega$, let 
$N = \bigcup\{N_\ell : \ell < \omega\}$. As this is a countable chain 
(recalling $\gK_\gs$ has amalgamation),
there is $q \in \clS(N)$ such that 
$\ell < \omega \Rightarrow q \rest N_\ell = q$ and as $q$ is
not algebraic (because each $q_n$ is not), and $\gs$ is full, clearly 
$q \in \clS_\gs(N)$; but $q$ contradicts the finite character of non-forking.  So for
some $n \ge 0$ we are stuck, but this cannot occur if $q_n \in \clS^\na_\gs(N_n)$.

[Why?  Because we are assuming that $\gs$ is type-full. Alternatively, we can use $\gs^\reg$, recalling that by \ref{a5}, we know that $\gs^\reg$ has enough regulars
and then we can apply \cite[8.3=L6.1tex]{Sh:705}.]

So for some $b \in N_n$ we have $q_n = \ortp(b,N_n,N_n)$;
i.e., $b$ realizes $q_n$ hence it realizes $p$.

Let $\bfJ = \{a_\ell:\ortp(a_\ell,N_\ell,N_{\ell +1})$ does not
fork over $N_0\}$.  By \cite[6.2]{Sh:705} we have 
$(M,N_n,\bfJ) = (N_0,N_n,\bfJ) \in K^{3,\vq}_\gs$
hence $\in K^{3,\qr}_\gs$ by \cite{Sh:705} so we are done.

\noindent 
2) Let $N,b,\bfJ$ be as in part (1) with $|\bfJ|$ minimal.  We
can find $N' \le_\gs N$ such that $(M,N',b) \in K^{3,\pr}_\gs$ and we 
can find $\bfJ'$ such that $\bfJ' \subseteq N'$ is 
independent regular in $(M,N')$ and maximal under those demands.  
Then we can find $N'' \le_\gs N'$ such that $(M,N'',\bfJ') \in
K^{3,\qr}_\gs$.  If $\ortp_\gs(b,N'',N') \in \clS^\na_\gs(N'')$ 
is not orthogonal to
$M$ we can contradict the maximality of $\bfJ'$ in $N'$ as in the
proof of part (1), 
so $\ortp_\gs(b,N'',N') \perp M$ (or $\notin \clS^\na_\gs(N)$).
Also without loss of generality $(N'',N',b) \in
K^{3,\pr}_\gs$, so by \cite{Sh:705} we 
have $(M,N',\bfJ') \in K^{3,\qr}_\gs$.
Hence there is an isomorphism $f$ from $N'$ onto $N''$ which is the
identity of $M \cup \bfJ'$ (by the uniqueness for
$K^{3,\qr}_\gs$).  So using $(N',f(b),\bfJ')$ for
$(N,b,\bfJ)$ we are done.

\noindent
3) If not, we can find $N_1,N_2,\bfJ_1,\bfJ_2,b$ such that $M
\le_\gs N_\ell \le_\gs N$ and the quadruple $(M,N_\ell,\bfJ_\ell,b)$ is
as in (a)-(e)+(f) of part (1)+(2) for $\ell =1,2$.  Assume toward
contradiction that $|\bfJ_1| \ne |\bfJ_2|$, so without loss of generality $|\bfJ_1| < |\bfJ_2|$. 

By ``$(M,N_\ell,b) \in K^{3,\pr}_\gs$," without loss of generality, $N_2 \le_\gs N_1$.

By \cite[10.15=L10b.11tex(3)]{Sh:705} for some 
$c \in \bfJ_2 \setminus \bfJ_1$, $\bfJ_1 \cup \{c\}$ is independent 
in $(M,N_1)$, in contradiction to $(M,N,\bfJ_1) \in K^{3,\vq}_\gs$ 
by \cite[10.15=L10b.11tex(4)]{Sh:705}.  

\noindent
4) Similarly.
\end{PROOF}

\begin{definition}\label{a11}
1) For $p \in \clS^\bs_\gs(M)$, let the \emph{weight} of $p$, $w(p)$,
be the unique natural number such that if $M \le_\gs M'$, $M'$ is 
brimmed, and $p' \in \clS^\bs_\gs(M')$ is a non-forking extension of $p$
\underline{then} it is the unique $|\bfJ|$ from Claim \ref{a8}(3). (It is 
a natural number.) 

\noindent
2) Let $w_\gs(a,M,N) = w(\ortp_\gs(a,M,N))$.
\end{definition}

\begin{claim}\label{w.3}
1) If $p \in \clS^\bs(M)$ is regular, \underline{then} $w(p) =1$. 

\noindent
2) If $\bfJ$ is independent in $(M,N)$ and $c \in N$, \underline{then} for
some $\bfJ' \subseteq \bfJ$ with $\le w_\gs(c,M,N)$
elements, $\{c\} \cup (\bfJ \setminus \bfJ')$ is independent in
$(M,N)$.
\end{claim}

\begin{PROOF}{\ref{w.3}}
Easy by now.
\end{PROOF}

\noindent
Note that the use of $\gC$ in Definition \ref{a17} is for
transparency only and can be avoided; see \ref{a26} below.

\begin{definition}\label{a17}
1) We say that $\bfP$ is an $M^*$-based family (inside $\gC$) \underline{when}:
\mn
\begin{enumerate}
    \item   $M^* <_{\gk[\gs]} \gC$ and $M^* \in K_\gs$
\sn
    \item  $\bfP \subseteq \bigcup\{\clS^\bs_\gs(M) : M \le_{\gk[\gs]} \gC$ and $M \in K_\gs\}$
\sn
    \item  $\bfP$ is preserved by automorphisms of $\gC$ over $M^*$.
\end{enumerate}
\mn
2) Let $p \in \clS^\bs_\gs(M)$, where $M \le_{\gk[\gs]}\gC$.
\mn
\begin{enumerate}
    \item We say that $p$ is hereditarily orthogonal to $\bfP$ (or $\bfP$-foreign) \underline{when}: 
    if $M \le_\gs N \le_{\gk[\gs]} \gC$, $q \in \clS^\bs_\gs(N)$, and $q \rest M = p$, \underline{then} $q$ is orthogonal to $p$.
\sn
    \item We say that $p$ is $\bfP$-regular \underline{when} $p$ is regular, not orthogonal to $\bfP$ and if $q \in \clS^\bs_\gs(M')$, $M \le_\gs M' <_{\gk[\gs]} \gC$, and $q$ is a forking extension of $p$ \underline{then} $q$ is hereditarily orthogonal to $\bfP$.
\sn
    \item  $p$ is weakly $\bfP$-regular \underline{if} it is regular and is not orthogonal to some $\bfP$-regular $p'$.
\end{enumerate}
\mn
3) $\bfP$ is normal \underline{when} $\bfP$ is a set of regular types
and each of them is $\bfP$-regular. 

\noindent
4) For $q \in \clS^\bs_\gs(M)$ and $M <_{\gk[\gs]} \gC$, 
let $w_\bfP(q)$ be defined as the natural number satisfying the following:
\mn
\begin{enumerate}
    \item[$\circledast$]  If $M \le_\gs M_1 \le_\gs M_2 \le_\gs \gC$, 
    $M_\ell$ is $(\lambda,*)$-brimmed, 
$b \in M_2$, $\ortp_\gs(b,M_1,M_2)$ is a non-forking extension 
of $q,(M_1,M_2,b) \in K^{3,\pr}_\gs$, $(M_1,M_2,\bfJ) \in 
K^{3,\qr}_\gs$, and $\bfJ$ is regular in 
$(M_1,M_2)$ (i.e. independent and $c \in \bfJ \Rightarrow \ortp_\gs
(c,M_1,M_2)$ is a regular type) \underline{then} 
$$w_\bfP(q) = \big|\{c \in \bfJ : \ortp_\gs(c,M_1,M) \text{ is weakly
$\bfP$-regular}\}\big|.$$
\end{enumerate}
\mn
5) We replace $\bfP$ by $p$ if $\bfP = \{p\}$, where $p \in
\clS^\bs(M^*)$ is regular (see \ref{a20}(1)).
\end{definition}

\begin{claim}\label{a20}
1) If $p \in \clS^\bs_\gs(M)$ is regular \underline{then}
$\{p\}$ is an $M$-based family and is normal. 

\noindent
2) Assume $\bfP$ is an $M^*$-based family.
If $q \in \clS^\bs_\gs(M)$ and $M^* \le_\gs M \le_{\gk[\gs]}\gC$ 
\underline{then} $w_\bfP(q)$ is well defined (and is a natural number). 

\noindent
3) In Definition \ref{a17}(4) we can find $\bfJ$ such that for every 
$c \in \bfJ_1$ we have: 
$$\ortp(c,M_1,M) \text{ is weakly $\bfP$-regular} \Rightarrow \ortp(c,M_1,M)\text{is $\bfP$-regular.}$$
\end{claim}

\begin{PROOF}{\ref{a20}}
 Should be clear. 
\end{PROOF}

\begin{discussion}\label{a23}  
1) It is tempting to try to generalize the notion of
$\bfP$-simple ($\bfP$-internal in Hrushovski's terminology)
and semi-regular.  An important property of those notions in the first
order case is that: e.g. 
\begin{itemize}
    \item[$(*)$] If $(\bar a/A) \not\perp p$ and $p$ regular, then for some equivalence relation $E$ definable over $A$, $\ortp(\bar a/E,A) \not\perp p$ and is $\{p\}$-simple.
\end{itemize}
The aim of defining $\{p\}$-simple is:
\begin{enumerate}
    \item For an element $a$ realizing $p$ over $A \subseteq \gC_T^\eq$, we can define the $p$-weight $w_p(a,A)$.

    \item The $p$-weight of such elements behaves like finite sequences from a vector space ---  so they behave like dimensions of vector spaces. 

    \item We have appropriate density results.
\end{enumerate}

2) First attempt towards (C) above: assume that $p,q \in \clS^\bs_\gs(M)$ are not
orthogonal, and we can define an equivalence relation ${\cE}^{p,q}_M$ on
$\{c \in \gC : c \text{ realizes } p\}$, defined by

\begin{equation*}
\begin{array}{clcr}
c_1\ {\cE}^{p,q}_M\ c_2 \quad \text{\underline{iff}} & \text{for every } d \in \gC 
\text{ realizing } q, \text{ we have} \\
  &\ortp_\gs(c_1 d,M,\gC) = \ortp_\gs (c_2d,M,\gC).
\end{array}
\end{equation*}

\mn
This (the desired property) may fail even in the first order case:  
suppose $p,q$ are definable over $a^* \in M$ (on
getting this, see later) and we have $\LL c_\ell : \ell \le n \RR$, 
$\LL M_\ell : \ell < n \RR$ such that 
$\ortp(c_\ell,M_\ell,\gC) = p_\ell$, each $p_\ell$
is parallel to $p$, $c_\ell\ {\cE}^{p,q}_{M_\ell}\ c_{\ell +1}$ but 
$c_0,c_n$ realize $p$ and $q$ respectively, and $\{c_0,c_n\}$ is 
independent over $M_0$.  Such a situation
defeats the attempt to define a $\bfP$-$\{q\}$-simple type $p/{\cE}$ as
in \cite[Ch.V]{Sh:c}.
 
However (see \cite[V, \S4]{Sh:c}), in first order logic we can find a saturated $N$ and  $a^* \in N$ such that
$$\ortp \big(M,\bigcup\limits_{\ell} M_\ell \cup\{c_0,\dotsc,c_n\} \big)$$ 
does not fork over $a^*$ and use ``average on the type with an ultrafilter
$c$ over $q(\gC) + a^*_t$" (for suitable $a^*_t$-s).  See more below.
\end{discussion}

\begin{discussion}\label{w20}
1) Assume ($\gs$ is full and) every $p \in \clS^\na_\gs(M)$ is 
representable by some $a_p \in M$ (e.g., in \cite{Sh:c}, the canonical
base $\Cb(p)$).  
We can define for $\bar a,\bar b \in {}^{\omega >} \gC$ 
when $\ortp(\bar a,\bar b,\gC)$ is stationary (and/or non-forking).  We 
should check the basic properties. See \S3.

\noindent
2) Assume $p \in \clS^\bs_\gs(M)$ is regular, 
definable over $\bar a^*$ (in the natural sense).  
We may wonder if the niceness of the 
dependence relation holds for $p \rest \bar a^*$?
\end{discussion}

If you feel that the use of a monster model is not natural in our
context, how do we ``translate" a set of types in $\gC^{\eq}$ 
preserved by every automorphism of $\gC$ which
is the identity on $A$? By using a ``place" defined by:

\begin{definition}\label{a26}
1) A \emph{local place} is a pair $\bfa = (M,A)$ such that $A \subseteq M
\in K_\gs$ (compare with \cite[\S1]{Sh:1239}).

\noindent
2) The places $(M_1,A_1),(M_2,A_2)$ are equivalent if $A_1 = A_2$ and
there are $n$ and $N_\ell \in K_\gs$ for $\ell \le n$
satisfying $A \subseteq N_\ell$
for $\ell = 0,\ldots,n$ such that $M_1 = N_0$, $M_2 = N_n$, and for each 
$\ell < n$, $N_\ell \le_\gs N_{\ell+1}$ or $N_{\ell+1} \le_\gs N_\ell$.  We write
$(M_1,A_1) \sim (M_2,A_2)$ or $M_1 \sim_{A_1} M_2$.

\noindent
3) For a local place $\bfa = (M,A)$, let $K_\bfa = K_{(M,A)} =
\{N : (N,A) \sim (M,A)\}$, so in $(M,A)/ {\sim}$ we fix both $A$ as a set
and the type it realizes in $M$ over $\varnothing$.

\noindent
4) We call such class $K_\bfa$ a \emph{place}.

\noindent
5) We say that $\bfP$ is an invariant set\footnote{Really a class.}
of types in a place $K_{(M,A)}$ \underline{when}:
\mn
\begin{enumerate}
    \item   $\bfP \subseteq \{\clS^\bs_\gs(N) : N \sim_A M\}$
\sn
    \item  Membership in $\bfP$ is preserved by isomorphism over $A$.
\sn
    \item If $N_1 \le_\gs N_2$ are both in $K_{(M,A)}$ and $p_2 \in \clS^\bs_\gs(N_2)$ does not fork over $N_1$ then $p_2 \in \bfP \Leftrightarrow p_2 \rest N_1 \in \bfP$. 
\end{enumerate}
\mn
6) We say $M \in K_\gs$ is brimmed over $A$ \underline{when} for some $N$ we
have $A \subseteq N \le_\gs M$ and $M$ is brimmed over $N$. 
\end{definition}

\begin{cd}\label{a32}
1) If $A \subseteq M \in K_\gs$ and $\bfP_0 \subseteq
   \clS^\bs_\gs(M)$ \underline{then} there is at most one invariant set
   $\bfP'$ of types in the place $K_{(M,A)}$ such that $\bfP^+
   \cap \clS^\bs_\gs(M) = \bfP_0$ and 
   $$M \le_\gs N \wedge p
   \in \bfP^+ \cap \clS^\bs_\gs(N) \Rightarrow ``p \text{ does not fork
   over } M".$$

\noindent
2) If in addition, $M$ is brimmed\footnote{$M$ is brimmed over $A$
  means that for some $M_1$, $A \subseteq M_1 \le_\gs M$ and $M$ is
  brimmed over $M_1$.}  over $A$ \underline{then} we can omit the
   last demand in part (1).

\noindent
3) If $\bfa = (M_1,A)$ and $(M_2,A) \in K_\bfa$ then $K_{(M_2,A)} = K_\bfa$. 
\end{cd}

\begin{PROOF}{\ref{a32}}
Easy.
\end{PROOF}

\begin{definition}\label{a35}
1) If in \ref{a32} there are such $\bfP$, we denote it by
   $\inv(\bfP_0) = \inv(\bfP_0,M)$.

\noindent
2) If $\bfP_0 = \{p\}$, then let $\inv(p) = \inv(p,M) = \inv(\{p\})$.

\noindent
3) We say $p \in \clS^\bs_\gs(M)$ \emph{does not split} (or \emph{is definable})
   over $A$ \underline{when} $\inv(p)$ is well defined.
\end{definition}

\newpage

\section {Imaginary elements, an essential-$(\mu,\lambda)$-AEC, and
frames} \label{2}

\subsection {Essentially AEC}\

We consider revising the definition of an AEC $\gk$, by allowing
function symbols in $\tau_\gk$ with infinite number of places
while retaining local characters, e.g., if $M_n \le M_{n+1}$ and 
$M = \bigcup\{M_n : n < \omega\}$ is uniquely determined.  Before this, we
introduce the relevant equivalence relations.  In this context, we can 
give name to equivalence classes for equivalence relations on
infinite sequences.

\begin{definition}\label{b2}
We say that $\gk$ is an essentially $[\lambda,\mu)$-AEC 
or ess-$[\lambda,\mu)$-AEC (or 
$[\lambda,\mu)$-EAEC\footnote{And we may write $(\mu,\lambda)$ instead of $[\lambda,\mu)$.})
\underline{if} ($\lambda < \mu$ and) it is an object consisting of:
\mn
\begin{enumerate}
    \item[\textbf{I.}]
    \begin{enumerate}
        \item A vocabulary $\tau = \tau_\gk$, which has predicates and function symbols of possibly infinite arity but $\le \lambda$.

        \item  A class $K = K_\gk$ of $\tau$-models.

        \item A two-place relation $\le_\gk$ on $K$.
    \end{enumerate}

    \sn 
    (Note that we allow $\mu = \infty$). 
\end{enumerate}
\mn
such that
\begin{enumerate}
    \item[\textbf{II.}]
    \begin{enumerate}
        \item If $M_1 \cong M_2$ \underline{then} 
        $M_1 \in K \Leftrightarrow M_2 \in K$.

        \item  If $(N_1,M_1) \cong (N_2,M_2)$ \underline{then} 
        $M_1 \le_\gk N_1 \Leftrightarrow M_2 \le_\gk N_2$.

        \item Every $M \in K$ has cardinality $\in [\lambda,\mu)$.

        \item  $\le_\gk$ is a partial order on $K$.
    \end{enumerate}

    \item[\textbf{III$_1$.}]  If $\LL M_i:i < \delta \RR$ is 
    $\le_\gk$-increasing and the cardinality of $\bigcup\limits_{i < \delta} M_i$ is less than $\mu$ 
    \underline{then} there is a unique $M \in K$ such that 
    $|M| = \bigcup\big\{|M_i| : i < \delta \big\}$ and 
    $i < \delta \Rightarrow M_i \le_\gk M$.

    \item[\textbf{III$_2$.}] If in addition, $i < \delta \Rightarrow M_i \le_\gk N$ \underline{then} $M \le_\gk N$.
    
    \item[\textbf{IV.}]  If $M_1 \subseteq M_2$ and $M_\ell \le_\gk N$ for 
    $\ell =1,2$ \underline{then} $M_1 \le_\gk M_2$.
    
    \item[\textbf{V.}]  If $A \subseteq N \in K$, \underline{then} there is $M$ 
    satisfying $A \subseteq M \le_\gk N$ and $\|M\| \le \lambda +|A|$. (Here it is enough to restrict ourselves to the case $|A| \le \lambda$.)
\end{enumerate}
\end{definition}

\begin{definition}
\label{b5}
1) We say $\gk$ is an ess-$\lambda$-AEC \underline{if} it is an
ess-$[\lambda,\lambda^+)$-AEC.

\noindent
2) We say $\gk$ is an ess-AEC \underline{if} there is $\lambda$ such
that it is an ess-$[\lambda,\infty)$-AEC, so $\lambda = 
\LST({\gk})$. 

\noindent
3) If $\gk$ is an ess-$[\lambda,\mu)$-AEC and 
$\lambda \le \lambda_1 < \mu_1 \le \mu$ then let 
$$K^\gk_{\lambda_1} = (K_\gk)_{\lambda_1} = \{M \in K_\gk : \|M\| = \lambda_1\}$$ 
and $K^\gk_{\lambda_1,\mu_1} = \{M \in K_\gk : \lambda_1 \le \|M\| < \mu_1\}$.

\noindent
4) We define $\Upsilon^{\oor}_\gk$ as in
   \cite[0.8=L11.1.3A(2)]{Sh:734}.

\noindent
5) We may omit the ``essentially" \underline{when} $\arity(\tau_\gk) =
\aleph_0$ where $\arity(\gk) = \arity(\tau_\gk)$ and for vocabulary 
$\tau$, 
$$\arity(\tau) = \min\{\kappa : \text{ every predicate and function
symbol has } \arity < \kappa\}.$$
\end{definition}

\noindent
We now consider the claims on ess-AECs.
\begin{claim}\label{b8}
Let $\gk$ be an ess-$[\lambda,\mu)$-AEC.  

\noindent
1) The parallel of $\Ax(III)_1,(III)_2$ holds with
a directed family $\LL M_t:t \in I \RR$.

\noindent
2) If $M \in K$ we can find $\LL M_{\bar a}:\bar a \in
{}^{\omega>}\!M \RR$ such that:
\mn
\begin{enumerate}
    \item   $\bar a \subseteq M_{\bar a} \le_\gk M$
\sn
    \item $\|M_{\bar a}\| = \lambda$
\sn
    \item If $\bar b$ is a permutation of $\bar a$ then $M_{\bar a} = M_{\bar b}$.
\sn
    \item  if $\bar a$ is a subsequence of $\bar b$ then $M_{\bar a} \le_\gk M_{\bar b}$.
\end{enumerate}
\mn
3) If $N \le_\gk M$ we can add in (2) that $\bar a \in
{}^{\omega >}\! N \Rightarrow M_{\bar a} \subseteq N$. 

\noindent
4) If for simplicity 
$$\lambda_* = \lambda + \sup\Big\{\textstyle\sum\big\{|R^M| : R \in \tau_\gk \big\} + \textstyle\sum \big\{|F^M| : F \in \tau_\gk \big\} : M \in K_\gk 
\text{ has cardinality }\lambda \Big\}$$ 
\underline{then} $K_\gk$ 
and $\{(M,N):N \le_\gk M\}$ are essentially $\PC_{\chi,\lambda_*}\!$-classes,
where $$\chi = \big|\{M/{\cong} : M \in K^\gk_\lambda\}\big|$$ noting that
$\chi \le 2^{2^\theta}$.  
That is, a sequence $\LL M_{\bar a} : \bar a\in {}^{\omega >}\! A\RR$
satisfying clauses (b),(c),(d) of part (2) such that 
$A = \bigcup \big\{|M_{\bar a}|:\bar a \in {}^{\omega >}\! A\big\}$ 
will represent a unique $M\in K_\gk$ with universe $A$ --- and similarly for
$\le_\gk$. (On the Definition of $\PC_{\chi,\lambda_*}$, see 
\cite[1.4(3)]{Sh:88r}.)  Note that if, in $\tau_\gk$, there are no two
distinct symbols with the same interpretation in every $M \in K_\gk$
then $|\tau){k_*}| \le 2^{2^\lambda}$.

\noindent
5) The results on omitting types in \cite{Sh:394} or
   \cite[0.9=L0n.8,0.2=0n.11]{Sh:734} hold, i.e., if
$\alpha < (2^{\lambda_*})^+ \Rightarrow
   K^\gk_{\beth_\alpha} \ne \varnothing$ \underline{then} $\theta \in
[\lambda,\mu) \Rightarrow K_\theta \ne \varnothing$ and there is an
$\EM$-model, i.e., $\Phi \in \Upsilon^{\oor}_\gk$ with $|\tau_\Phi| =
   |\tau_\gk| + \lambda$ and $\EM(I,\Phi)$ having cardinality
   $\lambda + |I|$ for any linear order $I$.

\noindent
6) The lemma on the equivalence of being
universal model homogeneous and of being saturated 
(see \cite[3.18=3.10]{Sh:300b} or \cite[1.14=L0.19]{Sh:600}) holds. 

\noindent
7) We can generalize the results of \cite[\S1]{Sh:600} on deriving an
ess-$(\infty,\lambda)$-AEC from an ess-$\lambda$-AEC.
\end{claim}

\begin{PROOF}{\ref{b8}}
The same proofs. On the generalization in \ref{b8}(7), see \cite[\S1]{Sh:1238}.
The point is that, in the language of \cite[\S1]{Sh:1238}, our $\gk$ is a
$(\lambda,\mu,\kappa)$-AEC (automatically with primes).
\end{PROOF}

\begin{remark}\label{b11}
1) In \ref{b8}(4) we can decrease the bound on $\chi$ if we have more nice
definitions of $K_\lambda$; e.g., if $\arity(\tau) \le \kappa$ then 
$\chi = 2^{(\lambda^{< \kappa}+|\tau|)}$, where $\arity(\tau) = \min\{\kappa :
\text{ every predicate and function symbol of $\tau$ has } \arity < \kappa\}$.

\noindent
2) We may use above $|\tau_\gs| \le \lambda$, $\arity(\tau_\gk) = \aleph_0$ to get that 
$$\big\{(M,\bar a)/{\cong} : M \in K^\gk_\lambda,\ \bar a \in {}^\lambda\! M \text{ lists } M\big\}$$
has cardinality $\le 2^\lambda$.  See also \ref{b56}.

\noindent
3) In \ref{b32} below, if we omit ``$\bbE$ is small" and $\lambda_1 =
   \sup\!\big\{|\seq(M)/\bbE_M|:M \in K^\gk_\lambda\big\}$ is $< \mu$ then
   $\gk_{[\lambda_1,\mu)}$ is an ess-$[\lambda_1,\mu)$-AEC.

\noindent
4) In Definition \ref{b2}, we may omit axiom V and define $\LST(\gk) \in
   [\lambda,\mu]$ naturally, and if $M \in K^\gk_\lambda \Rightarrow
   \mu > |\seq(M)/\bbE_M|$ \underline{then} in \ref{b32}(1) below we can omit
   ``$\bbE$ is small".

\noindent
5) Can we preserve in such ``transformation" the arity finiteness?  A 
natural candidate is trying to code $p \in \clS^\bs_\gs(M)$ 
by $\{\bar a : \bar a \in {}^{\omega >}\! M\}$, where
there are $M_0 \le_\gs M_1$ such that $M \le_\gs M_1$,
$\ortp(a_\ell,M_0,M_1)$ is parallel to $p$, and $\bar a$ is independent in
$(M_0,M_1)$.  If e.g., $K_\gs$ is saturated this helps but still
we suspect it may fail. 

\noindent
6) What is the meaning of ess-$[\lambda,\mu)$-AEC?
Can we look just at $\LL M_t : t \in I \RR$, $I$ directed, $t \le_I
s \Rightarrow M_t \le_\gs M_\gs \in K_\lambda$?  For
isomorphism types we take a kind of completion and so make more pairs
isomorphic, but $\bigcup\limits_{t \in I} M_t$ does not determine 
$\olsi M = \LL M_t : t \in I \RR$, and the completion may depend on this
representation.

\noindent
7) If we like to avoid this and this number is $\lambda'$, \underline{then}
 we should change the definition of $\seq(N)$ (see \ref{b14}(b)) to
 \begin{align*}
     \seq'(N) = \big\{\bar a : &\ \lh(\bar a) = \lambda,\ a_0 < \mu_*, \text{ and for some } M \le_\gs N \text{ from } K^\gk_\lambda, \\
     &\  \LL a_{1 + \alpha}:\alpha < \lambda\RR \text{ lists the members of $M$}  \big\}.
 \end{align*}

%$\seq'(N) = \big\{\bar a:\lh(\bar a) = \lambda$ and for some $M \le_\gs N$ from
%$K^\gk_\lambda$, $\LL a_{1 + \alpha}:\alpha < \lambda\RR$ lists
% the members of $M$ and $a_0 \in \{\gamma:\gamma <  \mu_*\} \big\}$. 
\end{remark}

\bigskip

\subsection {Imaginary Elements and Smooth Equivalence Relations}\
\bigskip

Now we return to our aim of getting canonical base for orbital types.

\begin{definition}\label{b14}
Let ${\gk} = (K_\gk,\le_\gk)$ be a $\lambda$-AEC or just
ess-$[\lambda,\mu)$-AEC (if $\gk_\lambda = \gk_\gs$ we may
write $\gs$ instead of $\gk_\lambda$, see \ref{b35}).  We say that 
$\bbE$ is a smooth $\gk_\lambda$-equivalence relation \underline{when}:
\mn
\begin{enumerate}
    \item   $\bbE$ is a function with domain $K_\gk$ mapping $M$ to $\bbE_M$.
\sn
    \item  For $M \in K_\gk$, $\bbE_M$ is an equivalence relation on a subset of 
    $$\seq(M) = \{\bar a \in {}^\lambda\! M : M \rest \Rang(\bar a) \le_\gk M\}$$ 
    so $\bar a$ is not necessarily without repetitions. Note that $\gk$ determines $\lambda$  (pedantically, when non-empty).
\sn
    \item  If $M_1 \le_\gk M_2$ then $\bbE_{M_2}
\rest \seq(M_1) = \bbE_{M_1}$.
\sn
    \item If $f$ is an isomorphism from $M_1 \in K_\gs$ onto $M_2$ \underline{then} $f$ maps $\bbE_{M_1}$ onto $\bbE_{M_2}$.

    \item If $\LL M_\alpha:\alpha \le \delta \RR$ is
$\le_\gs$-increasing continuous \underline{then} $$\big\{\bar a/\bbE_{M_\delta}:\bar a \in
\seq(M_\delta) \big\} = \big\{\bar a/\bbE_{M_\delta} : \bar a \in 
\textstyle\bigcup\limits_{\alpha < \delta} \seq(M_\alpha) \big\}.$$ 
\end{enumerate}
\mn
2) We say that $\bbE$ is small \underline{if} each $\bbE_M$ has $\le \|M\|$
 equivalence classes.
\end{definition}

\begin{remark}\label{b17}
1) Note that \underline{if} we have $\LL \bbE_i:i < i^* \RR$, where each
$\bbE_i$ is a smooth $\gk_\lambda$-equivalence relation and $i^* < \lambda^+$, 
\underline{then} we can find a smooth $\gk_\lambda$-equivalence relation $\bbE$ 
such that the $\bbE_M$-equivalence classes are essentially the $\bbE_i$-equivalence 
classes for $i < i^*$. In detail: without loss of generality $i^* \le \lambda$, and 
$\bar a\ \bbE_M\ \bar b$ \underline{iff} $\lh(\bar a) = \lh(\bar b)$ and
\begin{enumerate}
    \item[$\circledast_1$]   $i(\bar a) = i(\bar b)$, and if $i(\bar a) < i^*$ then
    $$\bar a \rest [1+i(\bar a)+1,\lambda)\ \bbE_{i(\bar a)}\ \bar b \rest [1+i(\bar b)+1,\lambda)$$ 
    where $i(\bar a) = \min\{j : j+1 < i(*) \wedge a_0 \ne a_{1+j} \text{ or } j = \lambda\}$.
\end{enumerate}

\mn
2) In fact $i^* \le 2^\lambda$ is OK: e.g., choose a function $\bfe$ 
from $\{e:e$ an equivalence relation on $\lambda\}$ to $i^*$. For
$\bar a,\bar b \in \seq(M)$ we let 
$i(\bar a) = \bfe\big(\{(i,j) : a_{2i+1} = a_{2j+1}\}\big)$ and
\mn
\begin{enumerate}
    \item[$\circledast_2$]   $\bar a\ \bbE_M\ \bar b$ \underline{iff}  $i(\bar a) = i(\bar b)$ and $\LL a_{2i} : i < \lambda \RR\ \bbE_{i(\bar a)}\ \LL b_{2i}:i < \lambda \RR$.
\end{enumerate}

\mn
3) We can redefine $\seq(M)$ as ${}^{\lambda \ge}\!M$, {but} then we have to make
minor changes above.
\end{remark}

\begin{definition}\label{b20}
Let $\gk$ be a $\lambda$-AEC or just ess-$[\lambda,\mu)$-AEC 
and $\bbE$ a small smooth $\gk$-equivalence relation and the reader
may assume for simplicity that the vocabulary $\tau_\gk$ has only predicates.
Also assume $F_*,c_*,P_* \notin \tau_\gk$.  
We define $\tau_*$ and $\gk_* = {\gk} \LL \bbE \RR =
(K_{\gk_*},\le_{\gk_*})$ as follows:
\mn
\begin{enumerate}
    \item  $\tau^* = \tau \cup \{F_*,c_*,P_*\}$ with $P_*$ a unary predicate, $c_*$ an individual constant and $F_*$ a $\lambda$-place function symbol.
\sn
    \item  $K_{\gk_*}$ is the class of $\tau^*_\gk$-models $M^*$ such that for some model 
    $M \in K_\gk$ we have:
    \begin{enumerate}
        \item[$(\alpha)$]   $|M| = P^{M^*}_*$
\sn
        \item[$(\beta)$]    If $R \in \tau$ then $R^{M^*} = R^M$.
\sn
        \item[$(\gamma)$]   If $F \in \tau$ has arity $\alpha$ then $F^{M^*} \rest M = F^M$ and for any $\bar a \in {}^\alpha(M^*)$, $\bar a \notin {}^\alpha\! M$ we have $F^{M^*}\!(\bar a) = c^{M^*}_*$ (or allow partial functions, or use $F^{M^*}\!(\bar a) = a_0$ when $\alpha > 0$ and $F^{M^*}(\LL\ \RR)$ when $\alpha = 0$; i.e. $F$ is an individual constant);
\sn
        \item[$(\delta)$]   $F_*$ is a $\lambda$-place function symbol, and:
        \begin{enumerate}
            \item[i.] If $\bar a \in \seq(M)$ then $F^{M^*}_*\!(\bar a) \in |M^*| \setminus |M| \setminus \{c^{M^*}_*\}$.
\sn
            \item[ii.] If $\bar a,\bar b \in \dom(\bbE) \subseteq \seq(M)$ then 
            $F^{M^*}_*\!(\bar a) = F^{M^*}_*\!(\bar b) \Leftrightarrow \bar a\ \bbE_M\ \bar b$.
\sn
            \item[iii.] If $\bar a \in {}^\lambda(M^*)$ and 
            $\bar a \notin \dom(\bbE) \subseteq \seq(M)$ then $F^{M^*}_*\!(\bar a) = c_*^{M^*}$.
        \end{enumerate}
\sn
        \item[$(\eps)$]   $c^{M^*}_* \notin |M|$, and if $b \in |M^*| \setminus |M| \setminus \{c^{M^*}_*\}$ then for some $\bar a \in \dom(\bbE) \subseteq \seq(M)$ we have 
        $F^{M^*}_*\!(\bar a) = b$.
    \end{enumerate}
\sn
    \item  $\le_{\gk_*}$ is the two-place relation on $K_{\gk_*}$ defined by: $M^* \le_{\gk_*} N^*$ \underline{if}
    \begin{enumerate}
        \item[$(\alpha)$]   $M^* \subseteq N^*$ and
\sn
        \item[$(\beta)$]   for some $M,N \in \gk$ as in clause (B) we have $M \le_\gk N$.
    \end{enumerate}
\end{enumerate}
\end{definition}

\begin{definition}\label{b26}
1) In \ref{b20}(1) we call $M \in {\gk}$ a \emph{witness} for 
$M^* \in K_{\gk_*}$ if it is as in clause (B) above.

\noindent
2) We call $M \le_\gk N$ a witness for $M^* \le_{\gk^*_\lambda} N^*$ if
they are as in clause (C) above.
\end{definition}

\begin{discussion}\label{b29}
Up to now we have restricted ourselves to vocabularies with each
predicate and function symbol of finite arity, and this restriction seems very
reasonable.  Moreover, it seems \emph{a priori} that for a parallel to
superstable, it is quite undesirable to have infinite arity.  Still,
our desire to have imaginary elements (in particular, canonical basis
for types) forces us to accept them.  The price is that for a general class of
$\tau$-models the union of increasing chains of $\tau$-models is not a well defined
$\tau$-model; more accurately, we can show its existence but not
smoothness. \underline{However}, inside the class $\gk\LL\bbE\RR$ defined above, it will be.
\end{discussion}

\begin{claim}\label{b32}
1) If $\gk$ is a $[\lambda,\mu)$-AEC or just an 
ess-$[\lambda,\mu)$-AEC and $\bbE$ a small smooth $\gk$-equivalence 
relation \underline{then} $\gk \LL \bbE \RR$
is an ess-$[\lambda,\mu)$-AEC.

\noindent
2) If $\gk$ has amalgamation and $\bbE$ is a small 
$\gk$-equivalence class \underline{then} ${\gk} \LL \bbE \RR$ has the
amalgamation property.
\end{claim}

\begin{PROOF}{\ref{b32}}
The same proofs.  Left as an exercise to the reader.
\end{PROOF}

\bigskip

\subsection {Good Frames}\
\bigskip

\noindent
Now we return to good frames.

\begin{definition}\label{b35}
We say that $\gs$ is a good ess-$[\lambda,\mu)$-frame \underline{if}
Definition \cite[2.1=L1.1tex]{Sh:600} is satisfied, except that:
\mn
\begin{enumerate}
    \item[(a)]   in clause (A), $\gK_\gs = (K_\gs,\le_\gs)$, $\gk$ is an ess-$[\lambda,\mu)$-AEC, and $\gK[\gs]$ is an ess-$(\infty,\lambda)$-AEC. 
\sn
    \item[(b)]   $K_\gs$ has a superlimit model in $\chi$ in every $\chi \in [\lambda,\mu)$.
\sn
    \item[(c)]  $K^\gs_\lambda/{\cong}$ has cardinality $\le 2^\lambda$, for convenience.
\end{enumerate}
\end{definition}

\begin{discussion}
\label{b40}
We may consider other relatives as our choice and mostly have similar
results.  In particular:
\mn
\begin{enumerate}
    \item   We can demand less: as in \cite[\S2]{Sh:842} we may replace $\clS^\bs_\gs$ by a formal version of $\clS^\bs_\gs$.
\sn
    \item   We may demand goodness only for $\gs_\lambda$, i.e. $\clS_\gs$ and $\nonfork{}{}_\gs$ apply only to models in $K_\lambda^\gs$ (hence we can use their nice properties from \cite[2.1=L1.1tex]{Sh:600}) and
    %$\gs$ \mgrimes{is the} restriction \mgrimes{of} the class of models to $K^\gs_\lambda$ and has only the formal properties above, so 
    amalgamation and JEP are required only for models of cardinality $\lambda$.
\end{enumerate}
\end{discussion}

\begin{claim}\label{b44}
All the definitions and results in \cite{Sh:600},
\cite{Sh:705} and \S1 here work for good ess-$[\lambda,\mu)$-frames.
\end{claim}

\begin{PROOF}{\ref{b44}}
No problem.  
\end{PROOF}

\begin{definition}\label{b47}
If $\gs$ is a $[\lambda,\mu)$-frame (see \ref{a2}), or just an ess-$[\lambda,\mu)$-frame,
and $\bbE$ a small smooth $\gs$-equivalence relation \underline{then} 
let ${\gt} = {\gs}\LL \bbE \RR$ be defined by:
\mn
\begin{enumerate}
    \item   $\gk_\gt = \gk_\gs\LL \bbE \RR$
\sn
    \item  $\clS^\bs_\gt(M^*) = \big\{\ortp_{\gk_\gt}(a,M^*,N^*) : M^* \le_{\gk_\gt} N^*$, and if $M \le_\gk N$ witness $M^*,N^* \in \gk_\gt$ then $a \in N \setminus M$ and $\ortp_\gs(a,M,N) \in \clS^\bs_\gs(M)\big\}$
\sn
    \item  Non-forking {is defined} similarly.
\end{enumerate}
\end{definition}

\begin{remark}\label{b50}
We may add: if $\gs$ is\footnote{The reader may ignore this version.}
 an NF-frame we define ${\gt} = {\gs}\LL \bbE \RR$ as an 
NF-frame similarly, see \cite{Sh:705}.
\end{remark}

\begin{claim}\label{b53}
1) If $\gs$ is a good ess-$[\lambda,\mu)$-frame, $\bbE$ a small, 
smooth $\gs$-equivalence relation \underline{then} ${\gs}\LL \bbE
\RR$ is a good ess-$[\lambda,\mu)$-frame. 

\noindent
2) In part (1), for every $\kappa$, $\dot I(\kappa,K^{\gs \LL\bbE\RR}) =
\dot I(\kappa,K^\gs)$. 

\noindent
3) If $\gs$ has primes/regulars \underline{then} ${\gs} \LL\bbE \RR$ does as well.
\end{claim}

\begin{remark}
\label{b54}  % 2023-03-08 06:12 3}
We may add: if $\gs$ is an NF-frame \underline{then} so is ${\gs}\LL \bbE
\RR$, hence $({\gs}\LL \bbE \RR)^{\full}$ is a
full NF-frame; see \cite{Sh:705}. 
\end{remark}

\begin{PROOF}{\ref{b53}}
 Straightforward.  
\end{PROOF}

\noindent
Our aim is to change $\gs$ inessentially such that for every
$p \in \clS^\bs_\gs(M)$ there is a canonical base,
etc.  The following claim shows that in the context we have presented
this can be done.

\begin{claim}\label{b56}
\textbf{The imaginary elements Claim}  

Assume $\gs$ a good $\lambda$-frame or just a good
ess-$[\lambda,\mu)$-frame.
 
\noindent
1) If $M_* \in K_\gs$ and $p^* \in \clS^\bs_\gs(M_*)$,
 \underline{then}\footnote{Note that there may well be an automorphism of $M^*$
 which maps $p^*$ to some $p^{**} \in \clS^\bs_\gs(M^*)$ such that
 $p^{**} \ne p^*$.}  \,  
there is a small, smooth $\gk_\gs$-equivalence relation $\bbE =
\bbE_{\gs,M_*,p^*}$ and function $\bfF$ such that:
\mn
\begin{enumerate}
    \item[$(*)$]   If $M_* \le_\gs N$, $\bar a \in \seq(N)$ (so 
    $M \defeq N \rest \Rang(\bar a) \le_\gs N$), and $M \cong M_*$, \underline{then}
    \begin{enumerate}
        \item[$(\alpha)$]  $\bfF(N,\bar a)$ is well defined iff $\bar a \in \dom(\bbE_N)$. 
        {If this is the case,} then $\bfF(N,\bar a)$ belongs to $\clS^\bs_\gs(N)$.
\sn
        \item[$(\beta)$]  $S \subseteq \big\{(N,\bar a,p) : N \in K_\gs,\ \bar a \in \dom(\bbE_N) \big\}$ is the minimal class such that:
        \begin{enumerate}
            \item[$(i)$]  If $\bar a \in \seq(M_*)$ and $p$ does not fork over 
            $M_* \rest \Rang(\bar a)$ then $(M_*,\bar a,p) \in S$.
\sn
            \item[$(ii)$]  $S$ is closed under isomorphisms.
\sn
            \item[$(iii)$]  If $N_1 \le_\gs N_2$, $p_2 \in \clS^\bs_\gs(N_2)$ does not fork over $\bar a \in \seq(N_1)$ then $(N_2,\bar a,p_2) \in S \Leftrightarrow (N_1,\bar a,p_2 \rest N_1) \in S$.
\sn
            \item[$(iv)$] If $\bar a_1,\bar a_2 \in \seq(N)$ and $p \in \clS^\bs_\gs(N)$ does not fork over $N \rest \Rang(\bar a_\ell)$ for $\ell=1,2$ then $(N_2,\bar a_1,p) \in S \Leftrightarrow (N_2,\bar a_2,p) \in S$.
        \end{enumerate}   
\sn
        \item[$(\gamma)$]  $\bfF(N,\bar a)=p$ iff $(N,\bar a,p) \in S$; hence if $\bar a,\bar b \in \seq(N)$ then $\bar a\ \bbE_N\ \bar b$ iff $\bfF(\bar a,N) = \bfF(\bar b,N)$.
    \end{enumerate}
\end{enumerate}
\mn
2) There are unique small\footnote{For `small' we use stability in $\lambda$.}
smooth $\bbE$-equivalence relations $\bbE_\gs$ and a function $\bfF$ such that:
\mn
\begin{enumerate}
    \item[$(**)(\alpha)$]   $\bfF(N,\bar a)$ is well defined iff $N \in K_\gs$ 
    and $\bar a \in \seq(N)$.
\sn
    \item[$(\beta)$]   $\bfF(N,\bar a)$, when defined, belongs to $\clS^\bs_\gs(N)$.
\sn
    \item[$(\gamma)$]  If $N \in K_\gs$ and $p \in \clS^\bs_\gs(N)$ \underline{then} there is 
    $\bar a \in \seq(N)$ such that $\Rang(\bar a) = N$ and $\bfF(N,\bar a) = p$.
\sn
    \item[$(\delta)$]  If $\bar a \in \seq(M)$ and $M \le_\gs N$ \underline{then} $\bfF(N,\bar a)$ is (well defined and is) the non-forking extension of $\bfF(M,\bar a)$.
\sn
    \item[$(\eps)$]  If $\bar a_\ell \in \seq(N)$ and $\bfF(N,\bar a_\ell)$ 
    is well defined for $\ell=1,2$ \underline{then} 
    $\bar a_1\ \bbE_N\ \bar a_2 \Leftrightarrow \bfF(N,\bar a_1) = \bfF(N,\bar a_2)$.
\sn
    \item[$(\zeta)$]   $\bfF$ commutes with isomorphisms.
\end{enumerate}
\mn
3) For ${\gt} = {\gs} \LL \bbE \RR$, where $\bbE$ is as in part (2), 
and whenever  % 2023-03-08 06:08 
$M^* \in K_\gt$ as witnessed by $M \in K_\gs$, and $p^* \in \clS^\bs_\gt(M^*)$ is 
projected to $p \in \clS^\bs_\gs(M)$, we let $\bas(p^*) = \bas(p) 
= \bfF(\bar a,M^*)/\bbE$ whenever $\bfF(M,\bar a) = p$.  
\mn
\begin{enumerate}
    \item[$(\alpha)$] If $M_\ell$ witnesses that $M^*_\ell \in K_\gt$ for $\ell=1,2$ and $(M^*_1,M^*_2,a) \in K^{3,\bs}_\gt$ {then} $(M_1,M_2,a) \in K^{3,\bs}_\gs$ and $p^* = \ortp_\gt(a,M^*_1,M^*_2)$, $p = \ortp_\gs(a,M_1,M_2)$. 

    \item[$(\beta)$] If $M^*_\ell \le_\gs M^*$ and $p_\ell \in \clS^\bs_\gt(M^*_\ell)$ 
    then $p^*_1 \parallel p^*_2 \Leftrightarrow \bas(p^*_1) = \bas(p^*_2)$.
\sn
    \item[$(\gamma)$]  $p^* \in \clS^\bs_\gt(M^*)$ does not split over $\bas(p^*)$ (see Definition \ref{a35}(3) or \cite[\S2 end]{Sh:705}).
\end{enumerate}
\end{claim}

\begin{PROOF}{\ref{b56}}
1) Let $M^{**} \le_\gs M^*$ be of cardinality
$\lambda$ such that $p^*$ does not fork over $M^{**}$.  Let $\bar a^*
= \LL a_\alpha:\alpha < \lambda\RR$ list the elements of
$M^{**}$.

We say that $p_1 \in \clS^\bs_\gs(M_1)$ is a \emph{weak copy} of $p^*$ when there is a witness $(M_0,M_2,p_2,f)$, which means:
\mn
\begin{enumerate}
    \item[$\circledast_1$] 
    \begin{enumerate}
        \item $M_0 \le_\gs M_2$ and $M_1 \le_\gs M_2$.

        \item if $\|M_1\| = \lambda$ then $\|M_2\| = \lambda$.

        \item $ f$ is an isomorphism from $M^{**}$ onto $M_0$.

        \item $ p_2 \in \clS^\bs_\gs(M_2)$ is a non-forking extension of $p_1$.

        \item $p_2$ does not fork over $M_0$.

        \item $f(p^* \rest M^{**})$ is $p_2 \rest M_0$.
    \end{enumerate}   
\end{enumerate}
\mn
For $M_1 \in K^\gs_\lambda$, $p_1 \in \clS^\bs_\gs(M_1)$ a weak copy of $p^*$, 
we say that $\bar b$ explicates its being a weak copy \underline{when}, for some witness $(M_0,M_2,p_2,f)$ and $\bar c$,
\mn
\begin{enumerate}
    \item[$\circledast_2$]  
    \begin{enumerate}
        \item $\bar b = \LL b_\alpha : \alpha < \lambda \RR$ lists the elements of $M_1$.

        \item $\bar c = \LL c_\alpha : \alpha < \lambda\RR$ lists the elements of $M_2$.

        \item $\{\alpha : b_{2 \alpha} = b_{2 \alpha+1}\}$ codes the following sets:
        \begin{enumerate}
            \item[$(\alpha)$]   The isomorphic type of $(M_2,\bar c)$.
\sn
            \item[$(\beta)$]    $\big\{(\alpha,\beta) : b_\alpha = c_\beta \big\}$
\sn
            \item[$(\gamma)$]    $\{(\alpha,\beta):f(a^*_\alpha) = c_\beta\}$
        \end{enumerate}
    \end{enumerate} 
\end{enumerate}
\sn
Now
\sn
\begin{enumerate}
    \item[$\circledast_3$]   If $p \in \clS^\bs_\gs(M)$ is a weak copy of $p^*$ 
    then for some $\bar a \in \seq(M)$, there is a $M_1 \le_\gs M$ over which $p$ does not
    fork such that $\bar a$ lists $M_1$ and explicates `$p \rest M_1$ is a weak copy of $p^*$.'
\sn
    \item[$\circledast_4$]   
    \begin{enumerate}
        \item If $M \in K^\gs_\lambda$ and $\bar b$ explicates `$p_1 \in \clS^\bs_\gs(M)$ is a weak copy of $p^*$,' \underline{then} we can reconstruct $p_1$ from $M$ and $\bar b$. (Call it $p_{M,\bar b}$.)

        %\mgrimes{Is $p^*$ a weak copy of itself? I think the first instance should be a $p_1$.}

        \item If $M \le_\gs N$, let $p_{N,\bar b}$ be its non-forking extension in $\clS^\bs_\gs(N)$. We also call it $\bfF(N,\bar b)$.
    \end{enumerate}
\end{enumerate}
\mn
Now we define $\bbE$. First, for $N \in K_\gs$ we define a
two-place relation $\bbE_N$.
\mn
\begin{enumerate}
    \item[$\circledast_5$]  
    \begin{enumerate}
        \item[$(\alpha)$] $\bbE_N$ is on $\{\bar a : \text{for some } M \le_\gs N$ of cardinality $\lambda$ and $p \in \clS^\bs_\gs(M)$ which is a copy of $p^*$, the sequence $\bar a$ explicates $p$ being  a weak copy of $p^*\}$.

        \item[$(\beta)$]   $\bar a_1\ \bbE_N\ \bar a_2$ iff $(\bar a_1,\bar a_2$ are as above and) $p_{N,\bar a_1} = p_{N,\bar a_2}$.
    \end{enumerate}
\end{enumerate}
\mn
Now
\mn
\begin{enumerate}
    \item[$\odot_1$]  For $N \in K_\gs$, $\bbE_N$ is an equivalence relation 
    on $\dom(E_N) \subseteq \seq(N)$.
\sn
    \item[$\odot_2$]  If $N_1 \le_\gs N_2$ and $\bar a \in \seq(N_1)$ 
    \underline{then} $\bar a \in \dom(\bbE_{N_1}) \Leftrightarrow \bar a \in \dom(\bbE_{N_2})$.
\sn
    \item[$\odot_3$]  If $N_1 \le_\gs N_2$ and $\bar a_1,\bar a_2 \in \dom(\bbE_{N_1})$ 
    then $\bar a_2\ \bbE_{N_1}\ \bar a_2 \Leftrightarrow \bar a_1\ \bbE_{N_2}\ \bar a_2$.
\sn
    \item[$\odot_4$]  If $\LL N_\alpha : \alpha \le \delta\RR$ is $\le_\gs$-increasing 
    continuous and $\bar a_1 \in \dom(\bbE_{N_\delta})$ then, for some $\alpha < \delta$ 
    and $\bar a_2 \in \dom(\bbE_{N_\alpha})$, we have $\bar a_1\ \bbE_{N_\delta}\ \bar a_2$.
\end{enumerate}
\mn
[Why?  Let $\bar a_2$ list the elements of $M_1 \le_\gs N_\delta$ and let 
$p = p_{N_\delta,\bar a_1}$ so $p \in \clS^\bs_\gs(N_\delta)$. Hence for some 
$\alpha < \delta$, $p$ does not fork over $M_\alpha$; hence for some 
$M'_1 \le_\gs M_\alpha$ of cardinality $\lambda$, the type $p$ does
not fork over $M'_1$.  Let $\bar a_2$ list the elements of $M'_1$ such
that it explicates $p \rest M'_1$ being a weak copy of $p^*$.  
So clearly $\bar a_2 \in \dom(\bbE_{N_\alpha}) \subseteq 
\dom(\bbE_{N_\delta})$ and $\bar a_1\ \bbE_{N_\delta}\ \bar a_2$.]

Clearly we are done.

\noindent
2) Similar, only we vary $(M^*,p^*)$ but it suffices to consider $2^\lambda$ such pairs.

\noindent
3) Should be clear.  
\end{PROOF}

\begin{dc}\label{b62}
Assume that $\gs$ is a good ess-$[\lambda,\mu)$-frame, so without loss of generality it is full.  
We can repeat the operations in \ref{b56}(3) and \ref{b53}(2), so 
after $\omega$ times we get ${\gt}_\omega$ which is full (that is, 
$\clS^\bs_{\gt_\omega}(M^\omega) = \clS^\na_{\gt_\omega}(M^\omega))$ 
and ${\gt}_\omega$ has canonical type-bases as witnessed by a 
function $\bas_{\gt_\omega}$ (see Definition \ref{b65}).
\end{dc}

\begin{PROOF}{\ref{b62}} % 2023-03-08 06:14 54}}
Should be clear.
\end{PROOF}

\begin{definition}\label{b65}
We say that $\gs$ has type bases \underline{if} there is a function
$\bas(-)$ such that:
\mn
\begin{enumerate}
    \item If $M \in K_\gs$ and $p \in \clS^\bs_\gs(M)$ \underline{then} $\bas(p)$ is (well defined and is) an element of $M$.
\sn
    \item $p$ does not split over $\bas(p)$; that, is any
    automorphism\footnote{There are reasonable stronger versions, but it follows that the function $\bas(-)$ satisfies them.} 
    of $M$ over $\bas(p)$ maps $p$ to itself.
\sn
    \item If $M \le_\gs N$ and $p \in \clS^\bs_\gs(N)$ then $\bas(p) \in M$ \underline{iff} $p$ does not fork over $M$.
\sn
    \item If $f$ is an isomorphism from $M_1 \in K_\gs$ onto $M_2 \in K_\gs$ 
    and $p_1 \in \clS^\bs(M_1)$ then $f(\bas(p_1)) = \bas(f(p_1))$.
\end{enumerate}
\end{definition}

\begin{remark}\label{b68}
1) In \S3 we can add:
\mn
\begin{enumerate}
    \item[(E)]   Strong uniqueness: if $A \subseteq M \le_{\gk(\gs)} \gC$ and $p \in \clS(A,\gC)$ is well defined \underline{then} $\bas(p) \in A$ and there is at most one $q \in \clS^\bs_\gs(M)$ such that $q$ extends $p$.  (Needed for non-forking extensions).
\end{enumerate}

2) In \ref{b71} we can work in $\gC$.
\end{remark}

\begin{definition}\label{b71}
We say that $\gs$ is equivalence-closed \underline{when}:
\mn
\begin{enumerate}
    \item  $\gs$ has type bases $p \mapsto \bas(p)$.
\sn
    \item  If $\bbE_M$ is a definition of an equivalence relation on ${}^{\omega >}\! M$ preserved by isomorphisms and $\le_\gs$-extensions (i.e. $M \le_\gs N \Rightarrow \bbE_M = \bbE_N \rest {}^{\omega >}\! M$) \underline{then} there is a definable function $F$ from ${}^{\omega >}\! M$ to $M$ such that $F^M(\bar a) = F^M(\bar b)$ iff $\bar a\ \bbE_M\ \bar b$. 
    %\saharon{or work in $\gC$}.
\end{enumerate}
\end{definition}

\noindent
To phrase the relation between $\gk$ and $\gk'$ we define the following. 

\begin{definition}\label{b77}
Assume $\gk_1,\gk_2$ are ess-$[\lambda,\mu)$-AECs.

\noindent
1) We say $\bfi$ is an interpretation in $\gk_2$ \underline{when}
$\bfi$ consists of
\mn
\begin{enumerate}
    \item A predicate $P^*_\bfi$.
\sn
    \item A subset $\tau_\bfi$ of $\tau_{\gk_2}$.
\end{enumerate}
\mn
2) In this case, for $M_2 \in K_{\gk_2}$, let $M^{[\bfi]}_2$ be the
$\tau_\bfi$-model $M_1 = M^{[\bfi]}_2$ with
\sn
\begin{itemize}
    \item universe $P^{M_2}_\bfi$, and
\sn
    \item  $R^{M_1} = R^{M_2} \rest |M_1|$ for $R \in \tau_\bfi$.
\sn
    \item  $F^{M_1}$ is defined similarly, so it can be a partial function even if $F^{M_2}$ is full.

    %\mgrimes{Something felt wrong here; just a hunch. Either ``$F^{M_{\red{1}}}$ can be a partial function even if $F^{M_2}$ is full,'' or some variant of ``$F^{M_{\red{\ell}}}$ can be a partial function even if $M_{\red{\ell}}$ is full'' would make perfect sense to me.}
\end{itemize}
\mn
3) We say that $\gk_1$ is $\bfi$-interpreted (or interpreted by
$\bfi$) in $\gk_2$ \underline{when}:
\mn
\begin{enumerate}
    \item   $\bfi$ is an interpretation in $\gk_1$.
\sn
    \item  $\tau_{\gk_1} = \tau_\bfi$
\sn
    \item  $K_{\gk_1} = \{M^{[\bfi]}_2:M_2 \in K_{\gk_2}\}$
\sn
    \item If $M_2 \le_{\gk_2} N_2$ then $M^{[\bfi]}_2 \le_{\gk_1} N^{[\bfi]}_2$.
\sn
    \item If $M_1 \le_{\gk_1} N_1$ and $N_1 = N^{[\bfi]}_2$ (so $N_2 \in K_{\gk_2}$) \underline{then} for some $M_2 \le_\gk N_2$ we have $M_1 = M^{[\bfi]}_2$.
\sn
    \item If $M_1 \le_{\gk_1} N_1$ and $M_1 = M^{[\bfi]}_2$ (so $M_2 \in K_{\gk_2}$) \underline{then} (possibly replacing $M_2$ by a model isomorphic to it over $M_1$) 
    there is $N_2 \in K_{\gk_2}$ such that 
    $M_2 \le_{\gk_2} N_2$ and $N_1 = N^{[\bfi]}_2$.
\end{enumerate}
\end{definition}

\begin{definition}\label{b80}
1) Assume $\gk_1$ is interpreted by $\bfi$ in $\gk_2$.  We say
   \emph{strictly} interpreted when: if $M^{[\bfi]}_2 = N^{[\bfi]}_2$
   then $M_2$ and $N_2$ are isomorphic over $M^{[\bfi]}_2$.

\noindent
2) We say $\gk_1$ is equivalent to $\gk_2$ if there are $n$
   and $\gk'_0,\dotsc,\gk'_n$ such that $\gk_1 = \gk'_0$, $\gk_2 =
   \gk'_n$ and for each $\ell < n$, $\gk_\ell$ is strictly interpreted in
   $\gk_{\ell +1}$ or vice versa.  Actually, we can demand $n=2$ and that
   $k_\ell$ is strictly interpreted in $\gk'_1$ for $\ell=1,2$.
\end{definition}

\begin{definition}\label{b83}
As above for (good) ess-$[\lambda,\mu)$-frames.
\end{definition}

\begin{claim}\label{b86}
Assume $\gs$ is a good $\ess$-$[\lambda,\mu)$-frame. 
\underline{Then} there exists $\gC$ (called a $\mu$-\emph{saturated monster} for $K_\gs$) such that: 
\mn
\begin{enumerate}
    \item[$(a)$]   $\gC$ is a $\tau_\gs$-model of cardinality $\le \mu$.
\sn
    \item[$(b)$]    $\gC$ is a union of some $\le_\gs$-increasing 
    continuous sequence $\LL M_\alpha : \alpha < \mu\RR$.
\sn
    \item[$(c)$]   if $M \in K_\gs$ (so $\lambda \le \|M\| < \mu$) \underline{then} 
    $M$ is $\le_\gs$-embeddable into some $M_\alpha$ from clause $(b)$.
\sn
    \item[$(d)$]    $M_{\alpha +1}$ is brimmed over $M_\alpha$ for $\alpha < \mu$.
\end{enumerate}
\end{claim}

\newpage

\section{$\bfP$-simple types} \label{3}

We define the basic types over sets not necessary models.  Note that
in Definition \ref{c11}(2) there is no real loss using $C$ of
cardinality $\in (\lambda,\mu)$, as we can replace $\lambda$ by
$\lambda_1 = \lambda + |C|$ and so replace $K_\gk$ to
$K^\gk_{[\lambda_1,\mu)}$.

\begin{hypothesis}\label{c2}
1) $\gs$ is a good ess-$[\lambda,\mu)$-frame (see Definition \ref{b35}).

\noindent
2) $\gs$ has type bases (see Definition \ref{b65}).

\noindent
3) $\gC$ will denote some $\mu$-saturated model for $K_\gs$ of
cardinality $\le \mu$; see \ref{b86}.

\noindent
4) But $M,A,\ldots$ will be $<_{\gk(\gs)} \gC$ and $\subseteq \gC$, respectively, but
 of cardinality $< \mu$.
\end{hypothesis}

\begin{definition}\label{c5}
Let $A \subseteq M \in K_\gs$.

\noindent
1) $\dcl(A,M) = \{a \in M : \text{if } M' \le_\gs M'',\ M \le_\gs M''$, 
and $A \subseteq M'$ then $a \in M'$ and for every automorphism $f$ 
of $M',\ f \rest A = \id_A \Rightarrow f(a) = a\}$.

\noindent
2) $\acl(A,M)$ is defined similarly, but only with the first demand.
\end{definition}

\begin{definition}\label{e.13y}
1) For $A \subseteq M \in K_\gs$ let
\[
\clS^\bs_\gs(A,M) = \{q \in \clS^\bs_\gs(M) : \bas(q) \in \dcl(A,\gC)\}.
\]

\sn
2) We call $p \in \clS^\bs_\gs(A,M)$ regular \underline{if} $p$ as a
member of $\clS^\bs_\gs(M)$ is regular. 
\end{definition}

\begin{definition}\label{c8}
1) $\bbE_\gs$ is as in Claim \ref{b56}(2).

\noindent
2) If $A \subseteq M \in K_\gs$ and $p \in \clS^\bs_\gs(M)$, 
then $p \in \clS^\bs_\gs(A,M)$ 
\underline{iff} $p$ is definable over $A$ (see \ref{a35}(3)) \underline{iff} 
$\inv(p)$ from Definition \ref{a35} is $\subseteq A$ and well defined.
\end{definition}

\begin{definition}\label{c11}
Let $A \subseteq \gC$. 

\noindent
1) We define a dependency relation on $\good(A,\gC) = \{c \in \gC: \text{for some }
M <_{\gk(\gs)} \gC,\ A \subseteq M$ and $\ortp(c,M,\gC)$ is
definable over some finite $\bar a \subseteq A\}$ as follows:
\mn
\begin{enumerate}
    \item[$\circledast$]    $c$ depends on $\bfJ$ in $(A,\gC)$ \underline{iff} there is no 
    $M <_{\gk(\gs)} \gC$ such that $A \cup \bfJ \subseteq M$ and $\ortp(c,M,\gC)$ 
    is the non-forking extension of $\ortp(c,\bar a,\gC)$, where $\bar a$ 
    witnesses $c \in \good(A,\gC)$. 
\end{enumerate}
\mn
2) We say that $C \in {}^{\mu >}[\gC]$ is \emph{good} over $(A,B)$
\underline{when} there is a brimmed $M <_{\gk(\gs)} \gC$ 
such that $B \cup A \subseteq M$ and $\ortp(C,M,\gC)$ (see
Definition \ref{a35}(3)) is definable over $A$. 
(In the first order context we could say $\{c,B\}$ is
independent over $A$, but here this is problematic as
$\ortp(B,A,\gC)$ is not necessarily basic.)

\noindent
3) We say $\LL A_\alpha:\alpha < \alpha^* \RR$ is independent
over $A$ in $\gC$ (see \cite[L8.8,6p.5(1)]{Sh:705}) \underline{if} we can 
find $M$ and $\LL M_\alpha : \alpha < \alpha^* \RR$ such that:
\mn
\begin{enumerate}
    \item[$\circledast$]
    \begin{enumerate}
        \item $A \subseteq M \le_{\gk(\gs)} M_\alpha <_\gs \gC$ for $\alpha < \alpha^*$.

        \item $M$ is brimmed.

        \item $A_\alpha \subseteq M_\alpha$

        \item $\ortp(A_\alpha,M,\gC)$ definable over $A$ (= does not split over $A$).

        \item $\LL M_\alpha : \alpha < \alpha^* \RR$ is independent over $M$.
    \end{enumerate} 
\end{enumerate}
\mn
3A) Similarly for ``over $(A,B)$". 

\noindent
4) We define `locally independent' naturally; that is, every finite subfamily is independent.
\end{definition}

\begin{claim}\label{c17}
Assume $a \in \gC$, $A \subseteq \gC$.

\noindent
1) $a \in \good(A,\gC)$ iff $a$ realizes $p \in \clS^\bs_\gs(M)$
   for some $M$ satisfying $A \subseteq M <_{\gk(\gs)} \gC$.
\end{claim}

\begin{claim}
\label{c20}
1) If $A_\alpha \subseteq \gC$ is good over 
$(A,\bigcup\limits_{i < \alpha} A_i)$ for $\alpha < \alpha^* < \omega$ 
then $\LL A_\alpha : \alpha < \alpha^* \RR$ is independent over $A$. 

\noindent
2) Independence is preserved by reordering. 

\noindent
3) If $p \in \clS^\bs_\gs(\bar a,\gC)$ is 
regular \underline{then} on $p(\gC) = \{c:c$ realizes $p\}$ the 
independence relation satisfies:
\mn
\begin{enumerate}
    \item[$(a)$]   Like clause (1).
\sn
    \item[$(b)$]   If $b^1_\ell$ depends on $\{b^0_0,\dotsc,b^0_{n-1}\}$ for 
    $\ell < k$ and $b^2$ depends on $\{b^1_\ell : \ell < k\}$ \underline{then} 
    $b^2$ depends on $\{b^0_\ell : \ell < n\}$.
\sn
    \item[$(c)$]  If $b$ depends on $\bfJ$, $\bfJ \subseteq \bfJ'$ \underline{then} $b$ depends on $\bfJ'$.
\end{enumerate}
\end{claim}

\begin{remark}\label{c23}
1) We have not mentioned finite character, but the local
independence satisfies it trivially.
\end{remark}

\begin{PROOF}{\ref{c20}}
Easy.
\end{PROOF}

\begin{definition}\label{c26}
1) Assume $q \in \clS^\bs_\gs(M)$ and $p \in \clS^\bs_\gs
(\bar a,\gC)$.  We say that $q$ is explicitly $(p,n)$-simple \underline{when}:
\mn
\begin{enumerate}
    \item[$\circledast$]   There are $b_0,\dotsc,b_{n-1},c$ such 
        that:\footnote{Clauses (c) + (e) are replacements for `$c$ is algebraic 
        over $\bar a + \{b_\ell:\ell < n\}$' and `each $b_\ell$ is necessary.'}
    \begin{enumerate}
        \item  $b_\ell$ realizes $p$.
\sn
        \item   $c$ realizes $q$.
\sn
        \item $b_\ell$ is not 
            good\footnote{`Not good' here is a replacement to 
            ``$\ortp(b_\ell,\bar a +c,\gC)$ does not fork over $\bar a$."} 
        over $(\bar a,c)$ for $\ell < n$.
\sn
        \item $\LL b_\ell:\ell < n \RR$ is independent over $\bar a$.
\sn
        \item $\LL c,b_0,\dotsc,b_{n-1} \RR$ is good over $\bar a$.
\sn
        \item If\footnote{This seems a reasonable choice here but we can take others; 
        this is an unreasonable choice for first order.}
        $c'$ realizes $q$ then $c=c'$ \underline{iff} for every $b \in p(\gC)$ 
        we have that $b$ is good over $(\bar a,c)$ iff $b$ is good over $(\bar a,c')$.
    \end{enumerate}
\end{enumerate}
\mn
1A) We say that $a$ is explicitly $(p,n)$-simple over $A$ \underline{if}
$\ortp(a,A,\gC)$ is; similarly, in the other definitions replacing
$(p,n)$ by $p$ will mean ``for some $n$." 

\noindent
2) Assume $q \in \clS^\bs_\gs(\bar a,\gC)$ and $\bfP$ as in Definition \ref{a17}.  
We say that $q$ is $\bfP$-simple if we can find $n$ and
explicitly $\bfP$-regular types $p_0,\dotsc,p_{n-1} \in \clS^\bs_\gs(\bar a,\gC)$ such that each $c \in p(\gC)$ 
is definable by its type over $\bar a \cup \bigcup\limits_{\ell <n} p_\ell(\gC)$.

\mn
3A) In part (1) we say weakly $(p,n)$-simple if in $\circledast$,
clause (f) is replaced by
\mn
\begin{enumerate}
    \item[(f)$'$]   If $b$ is good over $(\bar a,a^*_m)$ then $c$ and $c'$ 
    realize the same type over $\bar a \caret \LL a^*_m, b\RR$.
\end{enumerate}

\mn
3B) In part (1) we say $(p,n)$-simple \underline{if} for some 
$\bar a^* \in {}^{\omega >} \gC$ good over $\bar a$, for every 
$c \in q(\gC)$, there are $b_0,\dotsc,b_{n-1} \in p(\gC)$ such that 
$c \in \dcl(\bar a,\bar a^*,b_0,\dotsc,b_{n-1})$ and 
$\bar a \caret \LL b_0,\dotsc,b_{n-1}\RR$ is good over $\bar a$ if simple.

\noindent
4) Similarly in (2). 

\noindent
5) We define $\genw_p(b,\bar a)$ for $p$ regular and parallel to some $p' \in
\clS^\bs_\gs(\bar a)$. (Here $\genw$ stands for `general weight.')  Similarly
for $\genw_p(q)$.
\end{definition}

\noindent
We first list some obvious properties.

\begin{claim}\label{c32}
1) If $c$ is $\bfP$-simple over $\bar a$, with $\bar a \subseteq A 
\subset \gC$, \underline{then} $w_p(c,A)$ is finite. 

\noindent
2) The obvious implications.
\end{claim}

\begin{claim}\label{c35}
1) $[$Closures of the simple $\bs]$. 

\noindent
2) Assume $p \in \clS^\bs_\gs(\bar a,\gC)$.  
If $\bar b_1,\bar b_2$ are $p$-simple over $A$ \underline{then}
\mn
\begin{enumerate}
    \item[$(a)$]   $\bar b_1 \caret \bar b_2$ is $p$-simple (of course, $\ortp_\gs(\bar b_2 \bar b_2,\bar a,\gC)$ is not necessary in $\clS^\bs_\gs(\bar a,\gC)$ even if $\ortp_\gs(\bar b_\ell,\bar a,\gC) \in \clS^\bs_\gs(\bar a,\gC)$ for $\ell = 1,2$).
\sn
    \item[$(b)$]   Also, $\ortp(\bar b_2,\bar a b_1,\gC)$ is $\bfP$-simple.
\end{enumerate}
\mn
2) If $\bar b_\alpha$ is $p$-simple over $\bar a$ for $\alpha <
\alpha^*$ and $\pi : \beta^* \to \alpha^*$ one to one and onto, \underline{then}
$$\sum\limits_{\alpha < \alpha^*} \genw_p \big(b_\alpha,\bar a_* \cup
\textstyle\bigcup\limits_{\ell < \alpha} b_\alpha\big) = \sum\limits_{\beta < \beta^*}
\genw \big(b_{\pi(\beta)}, \bar a \cup \textstyle\bigcup\limits_{i < \beta} \bar b_{\pi(i)}\big).$$
\end{claim}

\mn
The following definition comes from \cite[6.9(1)=Lg29]{Sh:1238}.

\begin{definition}\label{g29}
Assume $p_1,p_2 \in \clS^\bs(M)$.  We say $p_1,p_2$ are
weakly orthogonal (and denote it $p_1 {\underset \wk \perp} p_2$) when the following implication holds: if
$M_0 \le_\gs M_\ell \le_\gs M_3$, $(M_0,M_\ell,a_\ell) \in
K^{3,\pr}_\gs$ and $\ortp_\gs(a_\ell,M_0,M_\ell) =
p_\ell$ for $\ell=1,2$ \underline{then} $\ortp_\gs(a_2,M_1,M_3)$ does not
fork over $M_0$. (this is symmetric by \textbf{Ax.E}(f).)
\end{definition}

%\mgrimes{All three of the orphaned references were to let the reader know the definition of `weakly orthogonal.' If it's that important, it needs to be in this paper.}

%\mgrimes{For reference, clauses (2)-(5) defined vanilla orthogonality (literally $p \perp q$), independence of a sequence over a triple $(M_0,M_1,M_2)$, $K^{3,\bs}_\gs$, and $K^{3,\qr}_\gs$. Check to see if any of those are necessary; you do say ``${\perp} = {\underset \wk \perp}$'' on the second line of this paper (in \ref{a2}).}

\begin{claim}\label{c38}
$[\gs$ is equivalence-closed.$]$ 

\noindent
Assume that $p,q \in \clS^\bs(M)$ are not weakly orthogonal (see \ref{g29}).
\underline{Then} for some $\bar a \in {}^{\omega >}\! M$ we have: $p,q$ are
definable over $\bar a$ (works without being stationary) and for some
$\gk_\gs$-definable function $\bfF$, for each $c \in q(\gC)$, 
$\ortp_\gs(\bfF(c,\bar a),\bar a,\gC) \in \clS^\bs_\gs(\bar a,\gC)$ 
and is explicitly $(p,n)$-simple for some $n$. (If, e.g., $M$ is $(\lambda,*)$-brimmed then $n = w_p(q)$.) 
\end{claim}

\begin{PROOF}{\ref{c38}}
We can find $n$ and $c_1,b_0,\dotsc,b_{n-1} \in
\gC$ with $c$ realizing $q$, $b_\ell$ realizing $p$, $\{b_\ell,c\}$
not independent over $M$, and $n$ maximal.  
Choose $\bar a \in {}^{\omega >}\! M$ such that 
$$\ortp_\gs(\LL c,b_0,\dotsc,b_{n-1}\RR,M,\gC)$$ 
is definable over $\bar a$.  
Define $E_{\bar a}$, an equivalence relation on $q(\gC)$: 
$c_1\ E_{\bar a}\ c_2$ iff for every $b \in p\gC$, we have ``$b$ is good over
$(a,c_1)" \Rightarrow ``b$ is good over $(\bar a,c_2)$."
By ``$\gs$ is eq-closed," we are done.
\end{PROOF}

\begin{claim}\label{c41}
1) Assume $p,q \in \clS^\bs_\gs(M)$ are weakly orthogonal 
(see \ref{g29}) 
but not orthogonal.  \underline{Then} we can
find $\bar a \in {}^{\omega >}\! M$ over which $p,q$ are definable and
$r_1 \in \clS^\bs_\gs(\bar a,\gC)$ such that letting
$p_1 = p \rest \bar a$, $q_1 = q \rest \bar a$, $n \defeq w_p(q) \ge 1$ we have:\footnote{We can say more concerning simple types.}
\mn
\begin{enumerate}
    \item[$\circledast^n_{\bar a,p_1,q_1,r_2}$]
    \begin{enumerate}
        \item[$(a)$] $p_1,q_1,r_1 \in \clS^\bs_\gs(\bar a)$, $\bar a \in{}^{\omega >}\gC$.
\sn
        \item[$(b)$]   $p_1,q_1$ are weakly orthogonal. 
        %(see e.g. Definition \ref{g29}(1)).
\sn
        \item[$(c)$]  If $\{a^*_n : n < \omega\} \subseteq r_1(\gC)$ is independent over $\bar a$ and $c$ realizes $q$ \underline{then} for infinitely many $m < \omega$ there is $b \in p(\gC)$ such that $b$ is good over $(\bar a,a^*_n)$ but not over $(\bar a,a^*_n,c)$.
\sn
        \item[$(d)$]  In (c) there really are $n$ independent $b_0,\dotsc,b_{n-1}$ which are all good over $(\bar a,a_n^*)$ but not over $(\bar a,a_n^*,c)$ (but we cannot find $n+1$ such $b$-s.).
    \end{enumerate}
\end{enumerate}

\mn
2) If $\circledast^n_{\bar a,p_1,q_1,r_1}$ then (see Definition
\ref{c26}(3) for some definable function $\bfF$, if $c$ realizes
$q_1$, $c^* = F(c,\bar a)$ and $\ortp_{\gn}(c^*,\bar a,\gC)$ is
$(p_1,n)$-simple. 
\end{claim}

\noindent
See proof below.

\begin{claim}\label{c44}
0) Assuming $A \subseteq \gC$ and $a \in A$, we say $\ortp(a,A,\gC)$ is \emph{finitary} \underline{when} it is definable over $A \cup \{a_0,\ldots,a_{n-1}\}$ for some $n$, where each $a_\ell$ is in $\gC$ and is good over $A$ inside $\gC$.

\noindent
1) If $a \in \dcl\big(\bigcup\{A_i : i < \alpha\}\cup A,\gC\big)$, $\ortp(a,A,\gC)$ is finitary, 
%(see Definition \saharon{xxx}; i.e. gotten by finitely many elements good ones over (A)\cyan{)} 
and $\{A_i : i < \alpha\}$ is independent over $A$ \underline{then} for some finite
$u \subseteq \alpha$ we have 
$$a \in \dcl \big(\textstyle\bigcup\{A_i : i \in u\} \cup A,\gC \big).$$ 

\noindent
2) If $\ortp(b,\bar a,\gC)$ is $\bfP$-simple, \underline{then} it is finitary.

\noindent
3) If $\{A_i : i < \alpha\}$ is independent over $A$ and $a$ is finitary
over $A$ then for some finite $u \subseteq \alpha$ (even $|u| < \wg(c,A)$), 
$\{A_i : i \in \alpha \setminus u\}$ is independent over $A,A \cup \{c\}$. (Or use $(A',A''),(A',A'' \cup \{c\})$.)
\end{claim}

\begin{definition}
\label{c50}
1) $\dcl(A) = \{a : f(a) = a \text{ for every automorphism $f$ of } \gC\}$. 

\noindent
2) $\dcl_\fin(A) = \bigcup\{\dcl(B) : B \subseteq A \text{ finite}\}$. 

\noindent
3) $a$ is finitary over $A$ \underline{if} there are $n < \omega$ and
$c_0,\dotsc,c_{n-1} \in$ $\good(A)$ such that $a \in \dcl(A \cup
\{c_0,\dotsc,c_{n-1}\})$. %\saharon{or $\dcl_\fin$?}. 

\noindent
4) For such $A$, let $\wg(a,A)$ be $w(\tp(a,A,\gC))$ when well defined.

\noindent
5) Strongly simple implies simple.
\end{definition}

\begin{claim}\label{c53}
In Definition \ref{c26}(3), for
some $m,k < \omega$ large enough, for every $c \in q(\gC)$ there
are $b_0,\dotsc,b_{m-1} \in \bigcup\limits_{\ell < n} p_\ell(\gC)$ such that 
$$c \in \dcl\big(\bar a \cup \{a^*_\ell:\ell < k\} \cup 
\{b_\ell:\ell < m\}\big).$$
\end{claim}

\begin{PROOF}{\ref{c41}}
Let $M_1,M_2 \in K_{\gs(\brim)}$ be
such that $M \le_\gs M_1 \le_\gs M_2$, $M_1$ is $(\lambda,*)$-brimmed over $M$, 
$p_\ell \in \clS^\bs_\gs(M_\ell)$ a 
non-forking extension of $p$, $q_\ell \in \clS^\bs_\gs(M_\ell)$ is 
a non-forking extension of $q$, $c \in M_2$ realizes $q_1$, and 
$(M_1,M_2,c) \in K^{3,\pr}_{\gs(\brim)}$.  Let
$b_\ell \in p_1(M_2)$ for $\ell < n^* \defeq w_p(q)$ be such that
$\{b_\ell : \ell < n^*\}$ is independent in $(M_1,M_2)$; let 
$\bar a^* \in {}^{\omega >}(M_1)$ be such that 
$\ortp_\gs(\LL c,b_0,\dotsc,b_{n-1} \RR,M_1,M_2)$ is definable over 
$\bar a^*$ and $r = \ortp_\gs(\bar a^*,M_1,M_2)$, $r^+ = \ortp(\bar
a^* \caret \LL b_0,\dotsc,b_{n-1} \RR,M,M_2)$.

Let $\bar a \in {}^{\omega >}\! M$ be such that 
$\ortp_\gs(\bar a^*,\LL c,b_0,\dotsc,b_{n-1} \RR,M,M_2)$ is definable over
$\bar a$.  As $M_1$ is $(\lambda,*)$-saturated over $M$, there is
$\{\bar a^*_f:f < \omega\} \subseteq r(\gC)$ independent in
$(M,M_1)$. Moreover, letting $a^*_\omega = \bar a^*$, we have 
$\LL a^*_\alpha : \alpha \le \omega \RR$ is independent in $(M,M_1)$.  
Clearly $\ortp_\gs(c \bar a^*_n,M,M_2)$ does not depend on $n$ hence we can
find $\big\LL \LL b^\alpha_\ell:\ell < n \RR : \alpha \le \omega\big\RR$ such that 
$b^\alpha_\ell \in M_2$, $b^\omega_\ell = b_\ell$, and
$\{c \bar a^*_\alpha,b^\alpha_0 \ldots b^\alpha_{n-1} : \alpha \le
\omega\}$. (As usual, this is because the index set is independent in $(M_1,M_2)$.)

The rest should be clear.
\end{PROOF}

\begin{definition}\label{c56}
Assume $\bar a \in {}^{\omega >} \gC$, $n < \omega$, and $p,q,r \in \clS^\bs(M)$ are as in the
definition of $p$-simple$^{[-]}$ but $p$ and $q$ are weakly
orthogonal (see e.g. Definition \ref{g29}(1)). Let $p$ be a definable related
function such that for any $\bar a^\nu_\ell \in r(\gC)$, 
$\ell < k^*$, the independent mapping 
$c \mapsto  \big\{b \in q(\gC) : R \gC \models R(b,c,\bar a^*_\ell) \big\}$ 
is a one-to-one function
from $q(\gC)$ into $$\big\{\LL J_\ell : \ell < k^* \RR : J_\ell
\subseteq p(\gC) \text{ is closed under dependence and has $p$-weight } 
n^*\big\}.$$ 

\noindent
1) We can define $E = E_{p,q,r}$, a two-place relation over 
$r(\gC)$: $\bar a^*_1\ E\ \bar a^*_2$ \underline{iff} $\bar a_1,\bar a_2 \in
r(\gC)$ have the same projection common to $p(\gC)$ and $q(\gC)$. 

\noindent
2) Define [\red{a} / \cyan{the}] unitless group on $r/E$ and its action on $q(\gC)$.
\end{definition}

\begin{remark}
\label{e.22}
1) A major point: as $q$ is $p$-simple,
$w_p(-)$ acts ``nicely" on $p(\gC)$, so if $c_1,c_2,c_3 \in q(\gC)$ then 
$w_p(\LL c_1,c_2,c_3 \RR \bar a) \le 3n^*$.  This
enables us to define averages using a finite sequence in a quite
satisfying way.  Alternatively, look more at averages of independent
sets. 

\noindent
2) \textbf{Silly Groups}:  Concerning interpreting groups, note that in our
present context, for every definable set $P^M$ we can add the group of
finite subsets of $P^M$ with symmetric difference (as addition).

\noindent
3) The axiomatization above has prototype $\gs$ where 
$K_\gs = \{M : M$ a $\kappa$-saturated model of $T\}$, 
${\le_\gs} = {\prec \rest K_\gs}$, $\nonfork{}{}_\gs$ is non-forking, 
$T$ a stable first order theory with $\kappa(T) \le \cf(\kappa)$.
But we may prefer to formalize the pair $(\gt,\gs)$ with $\gs$ as
above, $K_\gt =$ models of $T$, ${\le_\gt} = {\prec \rest K_\gt}$, 
$\nonfork{}{}_\gt$ is non-forking.

From $\gs$ we can reconstruct a $\gt$ by closing 
$\gk_\gs$ under direct limits, but in interesting cases we end up
with a bigger $\gt$.
\end{remark}

\bibliographystyle{amsalpha}
\bibliography{shlhetal}

\end{document}